\documentclass[11pt,a4paper]{article} 
\usepackage{graphicx,amsmath,bm}
\usepackage{amssymb,amscd}
\usepackage{caption2,color}

\usepackage[dvipdfm,
           pdfstartview=FitH,
           CJKbookmarks=true,
           bookmarksnumbered=true,
           bookmarksopen=true,
          colorlinks=true, 
          colorlinks=cyan,
           pdfborder=001,   
           citecolor=blue%
           ]{hyperref}

\makeatletter

\@addtoreset{equation}{section} \makeatother

\newtheorem{lemma}{Lemma}[section]

\newtheorem{rem}{Remark}[section]
\newtheorem{theorem}{Theorem}[section]
\newtheorem{corollary}{Corollary}[section]
\newcommand{\dx}{\,\mathrm{d}x}
\newcommand{\ds}{\,\mathrm{d}s}
\newcommand{\n}{\nabla}
\newcommand{\p}{\partial}

\newcommand{\seminorm}[1]{\lvert#1\rvert}
\newcommand{\rn}{\ensuremath{\mathbb{R}^N}}
\DeclareMathAlphabet{\mathsfsl}{OT1}{cmss}{m}{sl}

\renewcommand{\vec}[1]{\mbox{\boldmath$#1$}}
\newcommand{\oo}{\ensuremath{\Omega}}
\newcommand{\G}{\Gamma}
\newcommand{\defmath}{{\,\stackrel{\mbox{\rm\tiny def}}{=}\,}}
\newcommand{\diff}{\,\mathrm{d}}

\newcommand{\mdiv}{\,\mathrm{div}\,}
\newcommand{\mcurl}{\,\mathrm{curl}\,}

\newcommand{\md}{\mathrm{D}}
\newcommand{\mg}{\mathcal{G}}
\newcommand{\mt}{\mathcal{T}}
\newcommand{\mq}{\mathcal{Q}}
\newcommand{\mpp}{\mathcal{P}}
\newcommand{\my}{\mathcal{Y}}
\newcommand{\mv}{\mathcal{V}}

\newcommand{\vyt}{\vec{y}_t}

\newcommand{\vv}{\vec{V}}
\newcommand{\vw}{\vec{w}}
\newcommand{\vvv}{\vec{v}}

\newcommand{\vy}{\vec{y}}
\newcommand{\vyd}{\vec{y}_d}

\newcommand{\vu}{\vec{u}}
\newcommand{\vf}{\vec{f}}
\newcommand{\vg}{\vec{g}}
\newcommand{\vphi}{\vec{\varphi}}

\newcommand{\vpsi}{\vec{\psi}}

\newcommand{\vn}{\vec{n}}
\begin{document}

\title{Shape Optimization for Navier--Stokes Flow\footnote{This work was
supported by the National Natural Science Fund of China under grant
numbers 10371096 and 10671153 for ZM Gao and YC Ma.}}

\author{Zhiming Gao\thanks{Corresponding author.  School of Science, Xi'an Jiaotong University, P.O.Box 1844, Xi'an, Shaanxi, P.R.China, 710049. E--mail:dtgaozm@gmail.com.}\qquad
 Yichen Ma\footnote{School of Science, Xi'an Jiaotong University, Shaanxi, P.R.China, 710049. E-mail:\,ycma@mail.xjtu.edu.cn.}
 \qquad Hongwei Zhuang\thanks{Engineering College of Armed Police Force, Shaanxi,\,P.R.China, 710086.}}
\date{}
 \maketitle
\noindent{{\textbf{Abstract.\;}} This paper is concerned with the
optimal shape design of the newtonian viscous incompressible fluids
driven by the stationary nonhomogeneous Navier--Stokes equations. We
use three approaches to derive the structures of shape gradients for
some given cost functionals. The first one is to use the Piola
transformation and derive the state derivative and its associated
adjoint state; the second one is to use the differentiability of a
minimax formulation involving a Lagrangian functional with a
function space parametrization technique; the last one is to employ
the differentiability of a minimax formulation with a function space
embedding technique. Finally we apply a gradient type algorithm to
our problem and numerical examples show that our theory is useful
for
practical purpose and the proposed algorithm is feasible. \\[8pt]
{{\textbf{Keywords.\;}}
shape optimization; shape derivative; gradient algorithm; minimax formulation; material derivative; Navier-Stokes equations.\\[8pt]
{{\textbf{AMS(2000) subject classifications.\;}}35B37, 35Q30, 49K35, 49K40.

\section{Introduction}
This paper deals with the optimal shape design for the stationary Navier--Stokes flow. This problem is of great practical
importance in the design
and control of many industrial devices such as aircraft wings,
cars, turbines, boats, and so on. The control variable is the
shape of the fluid domain, the object is to minimize some cost functionals that
may be given by the designer, and finally we can obtain the optimal
shapes by numerical computation.

Optimal shape design has received considerable attention already. Early works concerning on existence of
solutions and differentiability of the quantity (such as, state, cost functional, etc.)
 with respect to shape deformation occupied most of the 1980s (see \cite{ce81,delfour,Piron,zolesio}), the stabilization of structures using boundary variation technique has been fully addressed in \cite{delfour,Piron,zolesio}. For the optimal shape design for Stokes flow, many people are contributed
to it, such as O.Pironneau \cite{piron74}, J.Simon \cite{si90}, ZM Gao \emph{et.al.}\cite{gao06,gao06b}, and so on.

In this paper, in order to derive the
structures of shape gradients with respect to the shape of the
variable domain for some given cost functionals in shape optimization problems for Navier--Stokes flow, we suggest the following three approaches:
\begin{itemize}
    \item [(i)] use the Piola transformation and derive the state derivative with respect to the shape of the fluid domain and its associated adjoint state;
    \item [(ii)]  utilize the differentiability of a minimax formulation involving a Lagrangian functional with a function space parametrization technique;
    \item [(iii)] employ the differentiability of a minimax formulation involving a Lagrangian functional with a function space embedding technique;
\end{itemize}

In \cite{gao-robin}, we use the first approach to solve a shape
optimization problem governed by a Robin problem, and in
\cite{gao06b}, we derive the expression of shape gradients for
Stokes optimization problem by the first approach. In this paper, we
use this approach to study the optimal shape design for Navier--Stokes
flow with small regularity data.

As we all known, many shape optimization problems can be expressed as a minimax of
some suitable Lagrangian functional. Theorems on the differentiability of a saddle point
(i.e., a minimax) of such Lagrangian functional with respect to a
parameter provides very powerful tools to obtain {shape gradients} by
{function space parametrization} or {function space embedding}
without the usual study of the state derivative approach.

The function space parametrization technique and function space
embedding technique are advocated by M.C.Delfour and
J.-P.Zol\'{e}sio to solving poisson equation with Dirichlet and
Nuemann condition (see\cite{delfour}). In our paper \cite{gao0601,gao06}, we
apply them to solve a Robin problem and a shape optimization problem for Stokes flow, respectively. However, in this paper we extend them to study the optimal shape design for Navier--Stokes flow in despite of its lack of rigorous mathematical justification in case where the Lagragnian formulation is not convex. We shall show how this theorem allows, at least formally to bypass the study of the differentiability of the state and obtain the expression of shape gradients with respect to the shape of the
variable domain for some given cost functionals.

We will find that the three approaches lead to the same expressions of the shape gradients for our given cost functionals. Hence, even if the last two approaches lacks from a rigorous mathematical framework, they allow more flexible computations which can be very useful for practical purpose. On the numerical point of view, we give the implementation of our problem in two dimensional case at the end of this paper, and the numerical results show that the last two approaches provide big efficiency for the shape optimization problem.

This paper is organized as follows. In \autoref{sec2}, we briefly
recall the velocity method which is used for the characterization of
the deformation of the shape of the domain and give the definitions
of Eulerian derivative and shape derivative. We also give the
description of the shape optimization problem for Navier--Stokes
flow.

In \autoref{sec:stated}, we prove the existence of the weak Piola
material derivative, and give the description of the shape
derivative. After that, we express the shape gradients of the
typical cost functionals $J_i(\oo)$, ($i=1,2$) by introducing the
corresponding adjoint state systems.

Section \ref{sec4} is devoted to the computation of the shape
gradient of the Lagrangian functional due to a minimax principle
concerning the differentiability of the minimax formulation by
{function space parametrization} technique and function space
embedding technique.

Finally in the last section, we give a gradient type
algorithm with some numerical examples to prove that our theory
could be very useful for the practical purpose and the proposed algorithm is feasible.

\section{Preliminaries and statement of the problem}\label{sec2}
\subsection{Elements of the velocity method}
Domains $\oo$ don't belong to a vector space and this requires the
development of {shape calculus} to make sense of a ``derivative" or
a ``gradient". To realize it, there are about three types of
techniques: J.Hadamard \cite{ha07}'s normal variation method, the
{perturbation of the identity} method by J.Simon \cite{si80} and the
{velocity method} (see J.Cea\cite{ce81} and
J.-P.Zolesio\cite{delfour,zo79}). We will use the velocity method
which contains the others. In that purpose, we choose an open set
$D$ in $\rn$ with the boundary $\p D$ piecewise $C^k$, and a
velocity space $\vec V\in \mathrm{E}^k :=\{\vec V\in
C([0,\varepsilon];\mathcal{D}^k(\bar{D},\rn)): \vec V\cdot\vn_{\p
D}=0\;\mbox{on }\p D\}$, where $\varepsilon$ is a small positive
real number and $\mathcal{D}^k(\bar{D},\rn)$ denotes the space of
all $k-$times continuous differentiable functions with compact
support contained in $\rn$ . The velocity field
$$\vec V(t)(x)=\vec V(t,x), \qquad x\in D,\quad t\geq 0$$
belongs to $\mathcal{D}^k(\bar{D},\rn)$ for each $t$. It can generate
transformations
 $$T_t(\vec V)X=x(t,X),\quad t\geq 0,\quad X\in D$$
through the following dynamical system
\begin{equation}\label{dynamical}
  \left\{%
  \begin{array}{ll}
  \frac{\diff x}{\diff t}(t,X)=\vec V(t,x(t))\\[3pt]
  x(0,X)=X
  \end{array}%
  \right.
\end{equation}
with the initial value $X$ given. We denote the "transformed domain"
$T_t(\vec V)(\oo)$ by $\oo_t(\vec V)$ at $t\geq 0$, and also set $\Gamma_t:=T_t(\Gamma)$.

There exists an interval $I=[0,\delta)$, $0<\delta\leq\varepsilon,$ and a one-to-one map $T_t$ from $\bar{D}$ onto $\bar{D}$ such that
\begin{itemize}
    \item [(i)] $T_0=\mathrm{I};$
    \item [(ii)] $(t,x)\mapsto T_t(x)$ belongs to $C^1(I;C^k(D;D))$ with $T_t(\p D)=\p D$;
    \item[(iii)]$(t,x)\mapsto T_t^{-1}(x)$ belongs to $C(I;C^k(D;D))$.
\end{itemize}
Such transformation are well studied in \cite{delfour}.

 Furthermore, for sufficiently small $t>0,$ the Jacobian $J_t$ is
 strictly positive:
 \begin{equation}\label{jacobian}
   J_t(x):=\det\seminorm{\md T_t(x)}=\det\md T_t(x)>0,
 \end{equation}
where $\md T_t(x)$ denotes the Jacobian matrix of the transformation
$T_t$ evaluated at a point $x\in D$ associated with the velocity
field $\vec V$. We will also use the following notation: $\md
T_t^{-1}(x)$ is the inverse of the matrix $\md T_t(x)$ , ${}^*\md
T_t^{-1}(x)$ is the transpose of the matrix $\md T_t^{-1}(x)$. These
quantities also satisfy the following lemma.
\begin{lemma}\label{lem:a}
    For any $\vec V\in E^k$, $\md T_t$ and $J_t$ are invertible. Moreover, $\md T_t$, $\md T_t^{-1}$ are in $C^1([0,\varepsilon];C^{k-1}(\bar{D};\mathbb{R}^{N\times N}))$, and $J_t$, $J_t^{-1}$ are in $C^1([0,\varepsilon];C^{k-1}(\bar{D};\mathbb{R}))$
\end{lemma}

Now let $J(\oo)$ be a real valued functional associated with any regular domain $\oo$, we say that this functional has a {\bf Eulerian derivative} at
$\oo$ in the direction $\vec V$ if the limit\\[6pt]
\begin{equation*}
\lim_{t\searrow 0}\frac{J(\oo_t)-J(\oo)}{t}:=\diff J(\oo;\vec
V)
\end{equation*}
exists.

 Furthermore, if the map
 $$\vec V\mapsto\diff
J(\oo;\vec V):\;\mathrm{E}^k\rightarrow\mathbb{R}$$ is linear and
continuous, we say that $J$ is {\bf shape differentiable} at
$\oo$. In the distributional sense we have
\begin{equation}\label{pri:shaped}
    \diff J(\oo;\vec V)=\langle \n J,\vec V\rangle_{\mathcal{D}^k(\bar{D},\rn)'\times \mathcal{D}^k(\bar{D},\rn)}.
\end{equation}
 When $J$ has a Eulerian derivative, we say that $\n J$ is the {\bf shape gradient} of $J$
at $\oo$.

Before closing this subsection, we introduce the following
functional spaces which will be used throughout this paper:
\begin{eqnarray*}
  H(\mdiv,\oo):=\{\vu\in L^2(\oo)^N: \;\mdiv\vu=0\mbox{ in
  }\oo,\;\vu\cdot\vn=0\mbox{ on }\p\oo\},\\
  H^1(\mdiv,\oo):=\{\vu\in H^1(\oo)^N:\;\mdiv\vu=0\mbox{ in
  }\oo\},\\
  H^1_0(\mdiv,\oo):=\{\vu\in H^1(\oo)^N:\;\mdiv\vu=0\mbox{ in
  }\oo,\;\vu|_{\p\oo}=0\}.
\end{eqnarray*}
\subsection{Statement of the shape optimization problem}
Let $\oo$ be the fluid domain in \rn ($N=2\;\mbox{or}\;3$), and the
boundary $\G:=\p\oo$. The fluid is described by its velocity $\vy$
and pressure $p$ satisfying the stationary Navier--Stokes equations:
\begin{equation}\label{ns:nonhomo}
  \left\{%
  \begin{array}{ll}
-\alpha\Delta\vy+\md\vy\cdot\vy+\n p=\vf&\quad\mbox{in}\;\oo\\
\mdiv \vy=0&\quad\mbox{in}\;\oo\\
\vy=\vg&\quad\mbox{on}\;\G
  \end{array}%
  \right.
\end{equation}
where $\alpha$ stands for the inverse of the Reynolds number whenever the variables are appropriately nondimensionalized,
$\vf$ denotes the given body force per unit mass, and $\vg$ is the given velocity at the boundary $\G$.

For the existence and uniqueness of the solution of the nonhomogeneous Navier--Stokes system (\ref{ns:nonhomo}), we have the following results (see \cite{temam01}).
\begin{theorem}
        \label{thm:ns}
        We suppose that $\oo$ is of class $C^1$. For
        \begin{eqnarray}
        &   \vf\in [L^2(\rn)]^N\\
        &\vg\in H^{ \frac{5}{2}}(\mdiv,\rn):=\left\{\vg\in\left [H^{\frac{5}{2}}\left(\rn\right)\right]^N:\;\mdiv\vg=0\right\},
        \end{eqnarray}
 there exists at least one $\vy\in H^1(\mdiv,\oo)$ and a distribution $p\in L^2(\oo)$
 on $\oo$ such that \eqref{ns:nonhomo} holds. Moreover,
 if $\alpha$ is sufficiently large and $\vg$ in $H^1-$norm is sufficiently small,
 there exists a unique solution $(\vy,p)\in H^1(\mdiv,\oo)\times L^2_0(\oo)$ of
 \eqref{ns:nonhomo}. In addition, if $\oo$ is of class $C^2$, we
 have $(\vy, p)\in (H^1(\mdiv,\oo)\cap H^2(\oo)^N)\times H^1(\oo)$.
\end{theorem}
We are interested in solving the following minimization problem
\begin{equation}\label{ns:cost}
 \min_{\oo\in\mathcal{O}} J_1(\oo)=\frac{1}{2}\int_{\oo}\seminorm{\vy-\vyd}^2\dx,
\end{equation}
or
\begin{equation}\label{ns:cost2}
    \min_{\oo\in\mathcal{O}} J_2(\oo)=\frac{\alpha}{2}\int_{\oo}\seminorm{\mcurl\vy}^2\dx.
\end{equation}
An example of the admissible set ${\mathcal{O}}$ is:
$$\mathcal{O}:=\{\oo\subset\rn: \seminorm{\tilde\G}=1\},$$
where $\tilde\G$ is the domain inside the closed boundary $\G$ and
$\seminorm{\tilde\G}$ is its volume or area in 2D.

\section{State derivative approach}\label{sec:stated}
In this section, we use the Piola transformation to bypass the
divergence free condition and then derive a weak material derivative
by the weak implicit function theorem. Then we will derive the
structure of the shape gradients of the cost functionals by
introducing the adjoint state equations associated with the
corresponding cost functional.

\subsection{Piola material derivative}

In order to deal with the nonhomogeneous Dirichlet boundary condition on $\Gamma$, we take $\tilde \vy=\vy-\vg$, where $\tilde\vy$ satisfies the following homogeneous Navier--Stokes system
\begin{equation}\label{ns:homo}
  \left\{%
  \begin{array}{ll}
-\alpha\Delta\tilde\vy+\md\tilde\vy\cdot\tilde\vy+\md\tilde\vy\cdot\vg+\md\vg\cdot\tilde\vy+\n p=\vec F&\quad\mbox{in}\;\oo\\
\mdiv \tilde\vy=0&\quad\mbox{in}\;\oo\\
\tilde\vy=0&\quad\mbox{on}\;\G
  \end{array}%
  \right.
\end{equation}
with $\vec F:=\vf+\alpha\Delta\vg+\md\vg\cdot\vg\in [L^2(\rn)]^N.$

We say that the function $\tilde\vy\in H^1_0(\mdiv,\oo)$ is called a weak solution of problem (\ref{ns:homo}) if it satisfies
\begin{equation}
    \label{ns:weak}
    \langle e(\tilde\vy),\vw\rangle=0,\qquad \vw\in H^1_0(\mdiv,\oo)
\end{equation}
with
\begin{equation}
    \label{ns:a}
    \langle e(\tilde\vy),\vw\rangle:=\int_\oo(\alpha\md\tilde\vy:\md\vw+\md\tilde\vy\cdot\tilde\vy\cdot\vw+\md\tilde\vy\cdot\vg\cdot\vw+\md\vg\cdot\tilde\vy\cdot\vw-\vec F\cdot\vw)\dx.
\end{equation}
As we all known, the divergence free condition coming from the fact
that the fluid has an homogeneous density and evolves as an
incompressible flow is difficult to impose on the mathematical and
numerical point of view. Therefore in order to work with the
divergence free condition, we need to introduce the following lemma
(see \cite{bios}).
\begin{lemma}
    The Piola transform
    \begin{eqnarray*}
        \Psi_t:\quad H^1(\mdiv,\oo)&\mapsto& H^1(\mdiv,\oo_t)\\
                             \vphi&\mapsto& ((J_t)^{-1}\md T_t\cdot\vphi)\circ T_t^{-1}
    \end{eqnarray*}
is an isomorphism.
\end{lemma}

Now by the transformation $T_t$, we consider the solution
$\tilde\vy_t$ defined on $\oo_t$ of the perturbed weak formulation:
\begin{equation}\label{piola:oot}
\int_{\oo_t}(\alpha\md\tilde\vy:\md\vw_t+\md\tilde\vy\cdot\tilde\vy\cdot\vw_t+\md\tilde\vy\cdot\vg\cdot\vw_t+\md\vg\cdot\tilde\vy\cdot\vw_t-\vec F\cdot\vw_t)\dx=0,\quad\forall\vw_t\in H^1_0(\mdiv,\oo_t),
\end{equation}
and introduce $\tilde\vy^t=\Psi_t^{-1}(\tilde\vy_t), \vw^t=\Psi_t^{-1}(\vw_t)$ defined on $\oo$.

We replace $\tilde\vy, \vw_t$ by $\Psi_t(\tilde\vy^t), \Psi_t(\vw^t)$ in the weak formulation (\ref{piola:oot}):
\begin{multline*}
    \int_{\oo_t}\left[\alpha\md(\Psi_t(\tilde\vy^t)):\md(\Psi_t(\vw^t))+\md(\Psi_t(\tilde\vy^t))\cdot\Psi_t(\tilde\vy^t)\cdot\Psi_t(\vw^t)+\md(\Psi_t(\tilde\vy^t))\cdot\vg\cdot\Psi_t(\vw^t)\right.\\\left.+\md\vg\cdot\Psi_t(\tilde\vy^t)\cdot\Psi_t(\vw^t)-\vec F\cdot\Psi_t(\vw^t)\right]\dx=0,\;\forall \vw^t\in H^1_0(\mdiv, \oo).
\end{multline*}
By the transformation $T_t$ and the following identities:
\begin{eqnarray*}
    \md(T_t^{-1})\circ T_t&=&\md T_t^{-1};\\
    \md(\vphi\circ T_t^{-1})&=&(\md\vphi\cdot\md T_t^{-1})\circ T_t^{-1};\\
    (\md\vg)\circ T_t&=& \md(\vg\circ T_t)\cdot\md T_t^{-1},
\end{eqnarray*}
we use a back transport in $\oo$ and obtain the following weak
formulation
\begin{equation}
    \langle e(t,\tilde\vy^t),\vw^t\rangle=0,\qquad \forall\vw^t\in H^1_0(\mdiv,\oo)
\end{equation}
with the notations
\begin{multline}\label{main:a}
\langle e(t,\vvv),\vw\rangle:=\alpha\int_\oo \md(B(t)\vvv):[\md(B(t)\vw)\cdot A(t)]\dx\\
+\int_\oo \md(B(t)\vvv)\cdot\vvv\cdot(B(t)\vw)\dx
+\int_\oo [\md(B(t)\vvv)\cdot\md T_t^{-1}]\cdot(\vg\circ T_t)\cdot(\md T_t\,\vw)\dx\\
+\int_\oo \md(\vg\circ T_t)\cdot\vvv\cdot(B(t)\vw)\dx -\int_\oo
(\vec F\circ T_t)\cdot(\md T_t\vw)\dx,
\end{multline}
and
\begin{equation*}
    A(t):=J_t\md T_t^{-1}{ }^*\md T_t^{-1};\qquad B(t)\vec \tau:=J_t^{-1}\md T_t\cdot\vec \tau.
\end{equation*}
Now we are interested in the derivability of the mapping
\begin{equation*}
t \mapsto \tilde\vy^t=\Psi_t^{-1}(\tilde\vy_t):\; [0,\varepsilon]\mapsto H^1_0(\mdiv,\oo)\\
\end{equation*}
where $\varepsilon>0$ is sufficiently small and $\tilde\vy^t\in H^1_0(\mdiv,\oo)$ is the solution of the state equation
\begin{equation}
    \langle e(t,\vvv),\vw\rangle=0,\qquad \forall \vw\in H^1_0(\mdiv,\oo).
\end{equation}
In order to prove the differentiability of $\tilde\vy^t$ with
respect to $t$ in a neighborhood of $t=0$, there maybe two
approaches:
\begin{itemize}
    \item [(i)] analysis of the differential quotient: $\lim\limits_{t\rightarrow 0}(\tilde\vy^t-\tilde\vy)/t$;
    \item [(ii)] derivation of the local differentiability of the solution $\tilde\vy$ associated to the implicit equation (\ref{ns:weak}).
\end{itemize}
We use the second approach. However, we can not use the classical implicit theorem, since it requires strong differentiability results in $H^{-1}$ for our case. Then we introduce the following weak implicit function theorem (see\cite{zo79}).
\begin{theorem}\label{theorem:weak}
    Let $X$, $Y'$ be two Banach spaces, $I$ an open bounded set in $\mathbb{R}$, and consider the map
    \begin{equation*}
      (t,x)\mapsto e(t,x):\;  I\times X\mapsto Y'
    \end{equation*}
    If the following hypothesis hold:
    \begin{itemize}
        \item [(i)]$t\mapsto \langle e(t,x),y\rangle$ is continuously differentiable for any $y\in Y$ and $(t,x)\mapsto\langle \p_t e(t,x),y\rangle$ is continuous;
        \item [(ii)] there exists $u\in X$ such that $u\in C^{0,1}(I;X)$ and $e(t,u(t))=0$, $\forall t\in I$;
        \item [(iii)] $x\mapsto e(t,x)$ is differentiable and $(t,x)\mapsto \p_x e(t,x)$ is continuous;
        \item [(iv)] there exists $t_0\in I$ such that $\p_x e(t,x)|_{(t_0,x(t_0))}$ is an isomorphism from $X$ to $Y'$,
    \end{itemize}
   the mapping
    \begin{equation*}
        t\mapsto u(t):\; I \mapsto X
    \end{equation*}
    is differentiable at $t=t_0$ for the weak topology in $X$ and its weak derivative $\dot{u}(t)$ is the solution of
    \begin{equation*}
        \langle \p_x e(t_0,u(t_0))\cdot\dot{u}(t_0),y\rangle+\langle \p_t e(t_0,u(t_0)),y\rangle=0,\quad \forall y\in Y.
    \end{equation*}
\end{theorem}
Now we state the main theorem of this subsection concerning on the differentiability of $\tilde\vy^t$ with respect to $t$.
\begin{theorem}\label{theorem:piola}
    The weak {\bf Piola material derivative}
    $$\dot{\tilde\vy}^P:=\lim_{t\rightarrow 0}\frac{\tilde\vy^t-\tilde\vy}{t}=\lim_{t\rightarrow 0}\frac{\Psi_t^{-1}(\tilde\vy_t)-\tilde\vy}{t}$$ exists and is characterized by the following weak formulation:
    \begin{equation}
        \langle \p_{\vec v} e(0,\vec v)|_{\vec v=\tilde\vy}\cdot{\dot{\tilde\vy}}^P,\vw\rangle
        +\langle \p_t e(0,\tilde\vy),\vw\rangle=0,\quad\forall \vw\in H^1_0(\mdiv,\oo),
    \end{equation}
i.e.,
\begin{multline}\label{piola:weak}
    \alpha\int_\oo\md\dot{\tilde\vy}^P:\md\vw\dx+\int_\oo\md\dot{\tilde\vy}^P\cdot\tilde\vy\cdot\vw+\md{\tilde\vy}\cdot\dot{\tilde\vy}^P\cdot\vw+\md\dot{\tilde\vy}^P\cdot\vg\cdot\vw+\md\vg\cdot\dot{\tilde\vy}^P\cdot\vw\dx\\=
\alpha\int_\oo\md\tilde\vy:\md\vw \mdiv\vec V\dx
+\alpha\int_\oo \md\tilde\vy:[-\md(\mdiv\vec V\vw)+\md(\md\vec V\vw)-\md\vw\md\vec V]\dx\\
+\alpha\int_\oo \md\vw:[-\md(\mdiv\vec V\tilde\vy)+\md(\md\vec V\tilde\vy)-\md\tilde\vy\md\vec V]\dx\\
+\int_\oo[\md(\md\vec V\tilde\vy-\mdiv\vec V\tilde\vy)\cdot\tilde\vy\cdot\vw+\md\tilde\vy\cdot\tilde\vy\cdot(\md\vec V\vw-\mdiv\vec V\vw)]\dx\\
+\int_\oo[\md(\md\vec V\tilde\vy-\mdiv\vec V\tilde\vy)-\md\tilde\vy\md\vec V]\cdot\vg\cdot\vw\dx
+\int_\oo \md\tilde\vy\cdot(\md\vg\vec V+{ }^*\md\vec V\vg)\cdot\vw\dx\\
+\int_\oo [\md(\md\vg\vec V)-\md\vg\mdiv\vec V+{ }^*\md\vec V\md\vg]\cdot\tilde\vy\cdot\vw\dx\\
-\int_\oo [\md(\vf+\alpha\Delta\vg+\md\vg\cdot\vg)\vec V+{}^*\md\vec V(\vf+\alpha\Delta\vg+\md\vg\cdot\vg)]\cdot\vw\dx.
\end{multline}
where $\tilde\vy$ is the solution of the weak formulation
\eqref{ns:weak}.
\end{theorem}
\noindent{\bf Proof.}\;
In order to apply Theorem \ref{theorem:weak}, we need to verify the four hypothesis of Theorem \ref{theorem:weak} for the mapping
\begin{equation*}
(t,\vvv)\mapsto e(t,\vvv):\;[0,\varepsilon]\times
H^1_0(\mdiv,\oo)\mapsto H^1_0(\mdiv,\oo)'.
\end{equation*}
To begin with, since the flow map $T_t\in C^1([0,\varepsilon];C^2(\bar{D},\rn))$ and Lemma \ref{lem:a}, the mapping
\begin{equation*}
    t\mapsto  \langle e(t,\vvv),\vw\rangle:  [0,\varepsilon]\mapsto  \mathbb{R}
\end{equation*}
is $C^1$ for any $\vvv, \vw\in H^1_0(\mdiv,\oo)$. On the other hand,
since $\vec F\in [L^2(\rn)]^N$, the mapping $t\mapsto \vec F\circ
T_t$ is only weakly differentiable in the space $[H^{-1}(\rn)]^N$,
thus the mapping $t\mapsto e(t,\vvv)$ is weakly differentiable. We denote by
$\p_t e(t,\vvv)$ its weak derivative. Since we have the following
three identities,
\begin{eqnarray}
    \label{deri:a}\frac{\diff}{\diff t}\md T_t&=&(\md\vec V(t)\circ T_t)\md T_t;\\
    \label{deri:b}\frac{\diff}{\diff t}J_t&=&(\mdiv \vec V(t))\circ T_t\, J_t;\\
    \label{deri:c}\frac{\diff}{\diff t}(\vf\circ T_t)&=&(\md\vf\cdot\vec V(t))\circ T_t,
\end{eqnarray}
the weak derivative  $\p_t e(t,\vvv)$ can be expressed as follows,
\begin{multline}
\langle\p_t e(t,\vvv),\vw\rangle=\alpha\int_\oo \md(B'(t)\vvv):[\md(B(t)\vw)\cdot A(t)]\dx\\
+\alpha\int_\oo\md(B(t)\vvv):[\md(B'(t)\vw)\cdot A(t)+\md(B(t)\vw)\cdot A'(t)]\dx\\
+\int_\oo [\md(B'(t)\vvv)\cdot\vvv\cdot(B(t)\vw)+\md(B(t)\vvv)\cdot\vvv\cdot(B'(t)\vw)]\dx\\
+\int_\oo [\md(B'(t)\vvv)\md T_t^{-1}-\md(B(t)\vvv)\cdot(\md T_t^{-1}(\md \vec V(t)\circ T_t)\md T_t^{-1})]\cdot (\vg\circ T_t)\cdot (\md T_t\vw)\dx\\
+\int_\oo[\md(B(t)\vvv)\md T_t^{-1}]\cdot \{[(\md\vg\cdot \vec V(t))\circ T_t]\cdot(\md T_t\vw)+(\vg\circ T_t)\cdot [(\md\vec V(t)\circ T_t)\md T_t\vw]\}\dx\\
+\int_\oo [\md(( \md\vg\cdot\vec V(t))\circ T_t)\cdot\vvv\cdot(B(t)\vw)+\md(\vg\circ T_t)\cdot\vvv\cdot(B'(t)\vw)]\dx\\
-\int_\oo [{ }^*\md\vec V(t)\vec F+\md\vec F\vec V(t)]\circ T_t\md T_t\,\vw\dx,
\end{multline}
where the notation
\begin{eqnarray*}
    &B'(t)\vec\tau:=\frac{\p}{\p t}\left[B(t)\vec\tau\right]=[\md\vec V(t)\circ T_t-(\mdiv \vec V(t)\circ T_t)\mathrm{I}]B(t)\vec \tau;\\[6pt]
    &A'(t):=\frac{\p}{\p t}A(t)=[\mdiv \vec V(t)\circ T_t-\md T_t^{-1}\md\vec V(t)\circ T_t]\,A(t)-{}^*[\md T_t^{-1}\md\vec V(t)\circ T_t\, A(t)].
\end{eqnarray*}
It is easy to check that the mapping $(t,\vvv)\mapsto \p_t e(t,\vvv)$ is
weakly continuous from $[0,\varepsilon]\times H^1_0(\mdiv,\oo)$ to
$H^1_0(\mdiv,\oo)'$.

We set $t=0$, use $T_t|_{t=0}=\mathrm{I}$ and $\vec V(t)|_{t=0}=\vec V$, then obtain
\begin{multline}\label{zero}
    \langle \p_t e(0,\vvv),\vw\rangle\\=\alpha\int_\oo\md\vvv:\md\vw \mdiv\vec V\dx
+\alpha\int_\oo \md\vvv:[-\md(\mdiv\vec V\vw)+\md(\md\vec V\vw)-\md\vw\md\vec V]\dx\\
+\alpha\int_\oo \md\vw:[-\md(\mdiv\vec V\vvv)+\md(\md\vec V\vvv)-\md\vvv\md\vec V]\dx\\
+\int_\oo[\md(\md\vec V\vvv-\mdiv\vec V\vvv)\cdot\vvv\cdot\vw+\md\vvv\cdot\vvv\cdot(\md\vec V\vw-\mdiv\vec V\vw)]\dx\\
+\int_\oo[\md(\md\vec V\vvv-\mdiv\vec V\vvv)-\md\vvv\md\vec
V]\cdot\vg\cdot\vw\dx
+\int_\oo \md\vvv\cdot(\md\vg\vec V+{ }^*\md\vec V\vg)\cdot\vw\dx\\
+\int_\oo [\md(\md\vg\vec V)-\md\vg\mdiv\vec V+{ }^*\md\vec V\md\vg]\cdot\tilde\vy\cdot\vw\dx
-\int_\oo (\md\vec F\vec V+{}^*\md\vec V\vec F)\cdot\vw\dx.
\end{multline}
To verify (ii), we follow the same steps described in R.Dziri\cite{dziri95} to find an identity satisfied by $\tilde\vy^{t_1}-\tilde\vy^{t_2}$ and prove that the solution $\tilde\vy^t\in H^1_0(\mdiv,\oo)$ of the weak formulation
\begin{equation*}
\langle e(t,\vvv),\vw\rangle=0,\qquad \forall \vw\in
H^1_0(\mdiv,\oo)
\end{equation*}
is Lipschitz with respect to $t$.

It is easy to find that $\vvv\mapsto e(t,\vvv)$ is differentiable,
and the derivative of $e(t,\vvv)$ with respect to $\vvv$ in the
direction $\delta\vvv$ is
\begin{multline*}
\langle\p_{\vvv} e(t,\vvv)\cdot\delta\vvv,\vw\rangle=\int_\oo\alpha\md(B(t)\delta\vvv):[\md(B(t)\vw)A(t)]\dx\\
+\int_\oo[\md(B(t)\delta\vvv)\cdot\vvv\cdot(B(t)\vw)+\md(B(t)\vvv)\cdot\delta\vvv\cdot(B(t)\vw)]\dx\\
+\int_\oo\{[\md(B(t)\delta\vvv)\md T_t^{-1}]\cdot(\vg\circ T_t)\cdot
(\md T_t\vw)+\md(\vg\circ T_t)\cdot\delta\vvv\cdot(B(t)\vw)\}\dx.
\end{multline*}
Obviously, $\p_{\vvv} e(t,\vvv)$ is continuous, and when we take
$t=0$,
\begin{multline}\label{piola:dy}
    \langle\p_{\vvv} e(0,\vvv)\cdot\delta\vvv,\vw\rangle=\alpha\int_\oo\md\delta\vvv:\md\vw\dx+\int_\oo\md(\delta\vvv)\cdot\vvv\cdot\vw+\md\vvv\cdot\delta\vvv\cdot\vw\\+\md(\delta\vvv)\cdot\vg\cdot\vw+\md\vg\cdot\delta\vvv\cdot\vw\dx.
\end{multline}
Furthermore, the mapping $\delta\vvv\mapsto \p_{\vvv}
e(0,\vvv)\cdot\delta\vvv$ is an isomorphism from $H^1_0(\mdiv,\oo)$
to its dual. Indeed, this result follows from the uniqueness and
existence of the Navier--Stokes system, i.e., Theorem \ref{thm:ns}.

Finally, all the hypothesis are satisfied by (\ref{main:a}), we can apply Theorem \ref{theorem:weak} to (\ref{main:a})
and then use (\ref{zero}) and (\ref{piola:dy}) to obtain (\ref{piola:weak}).\hfill $\square$

\subsection{Shape derivative}
In this subsection, we will characterize the shape derivative $\tilde\vy'$, i.e., the derivative of the state $\tilde\vy$ with respect to the shape of the domain.
\begin{theorem}\label{ns:shaped}
    Assume that $\oo$ is of class $C^2$, $\tilde\vy\in H^1_0(\mdiv,\oo)$ solves
     the weak formulation (\ref{ns:weak}) in $\oo$ and $\tilde\vy_t\in H^1_0(\mdiv,\oo_t)$
     solves the perturbed weak formulation (\ref{piola:oot}) in $\oo_t$, then the {\bf shape derivative}
    $$\tilde\vy':=\lim_{t\rightarrow 0}\frac{\tilde\vy_t-\tilde\vy}{t}$$
    exists and is characterized as the solution of
\begin{equation}\label{ns:shape}
    \left\{
    \begin{array}{lll}
        -\alpha\Delta\tilde\vy'+\md\tilde\vy'\cdot\tilde\vy+\md\tilde\vy\cdot\tilde\vy'+\md\tilde\vy'\cdot\vg+\md\vg\cdot\tilde\vy'+\n p'=0 &\qquad\mbox{in }\oo\\
        \mdiv \tilde\vy'=0 &\qquad\mbox{in  }\oo\\
        \tilde\vy'=-(\md\tilde\vy\cdot\vn)\vv_n&\qquad\mbox{on  }\Gamma
    \end{array}
    \right.
\end{equation}
\end{theorem}
\noindent{\bf Proof.}\;  We recall that $\tilde\vy_t\in H^1_0(\mdiv,\oo_t)$ satisfies the following weak formulation
\begin{equation}\label{stokes:c}
    \int_{\oo_t}(\alpha\md\tilde\vy_t:\md\vw+\md\tilde\vy_t\cdot\tilde\vy_t\cdot\vw+\md\tilde\vy_t\cdot\vg\cdot\vw+\md\vg\cdot\tilde\vy_t\cdot\vw)\dx-\int_{\oo_t}\vec F\cdot\vw\dx=0
\end{equation}
for any $\vw\in H^1_0(\mdiv,\oo_t)$.

To begin with, we introduce the following Hadamard formula (see
\cite{delfour,zolesio})
\begin{equation}\label{hadamard}
 \frac{\diff{}}{\diff t}\int_{\oo_t}F(t,x)\dx=\int_{\oo_t}
 \frac{\p F}{\p t}(t,x)\dx+\int_{\p\oo_t} F(t,x)\,\vec
 V\cdot\vn_t\diff\G_t,
 \end{equation}
 for a sufficiently smooth functional
$F:[0,\tau]\times\rn\rightarrow\mathbb{R}$.

Now we set a function $\vphi\in\mathcal{D}(\oo)^N$ and
$\mdiv\vphi=0$ in $\oo$. Obviously when $t$ is sufficiently small,
$\vphi$ belongs to the sobolev space $H^1_0(\mdiv,\oo_t)$. Hence we
can use (\ref{hadamard}) to differentiate \eqref{stokes:c} with
$\vw=\vphi$,
\begin{multline*}
    \int_\oo(\alpha\md\tilde\vy':\md\vphi+\md\tilde\vy'\cdot\tilde\vy\cdot\vphi+\md\tilde\vy\cdot\tilde\vy'\cdot\vphi
    +\md\tilde\vy'\cdot\vg\cdot\vphi+\md\vg\cdot\tilde\vy'\cdot\vphi)\dx\\
    +\int_\Gamma (\alpha\md\tilde\vy:\md\vphi+\md\tilde\vy\cdot\tilde\vy\cdot\vphi+\md\tilde\vy\cdot\vg\cdot\vphi
    +\md\vg\cdot\tilde\vy\cdot\vphi-\vec F\cdot\vphi)\vv_n\diff s=0.
\end{multline*}
Since $\vphi$ has a compact support, the boundary integral vanishes.
Using integration by parts for the first term in the distributed
integral, we obtain
\begin{equation}
    \int_\oo(-\alpha\Delta\tilde\vy'+\md\tilde\vy'\cdot\tilde\vy+\md\tilde\vy\cdot\tilde\vy'+\md\tilde\vy'\cdot\vg
    +\md\vg\cdot\tilde\vy')\cdot\vphi\dx=0.
\end{equation}
Then there exists some distribution $p'\in L^2(\oo)$ such that
$$-\alpha\Delta\tilde\vy'+\md\tilde\vy'\cdot\tilde\vy
+\md\tilde\vy\cdot\tilde\vy'+\md\tilde\vy'\cdot\vg+\md\vg\cdot\tilde\vy'=-\n p'$$
in the distributional sense in $\oo$.

Now we recall that for each $t$, $\Psi_t^{-1}(\tilde\vy_t)$ belongs
to the Sobolev space $H^1_0(\mdiv,\oo)$, then we can deduce that its material
derivative vanishes on the boundary $\G$. Thus we obtain the
shape derivative of $\tilde\vy$ at the boundary,
\begin{equation*}
    \tilde\vy'=-\md\tilde\vy\cdot\vv,\qquad\mbox{on  }\Gamma
\end{equation*}
Since $\tilde\vy|_\G=0,$ we have $\md\tilde\vy|_\G=\md\tilde\vy\cdot\vn^*\vn,$ and then
\begin{equation*}
    \tilde\vy'=-(\md\tilde\vy\cdot\vn)\vv_n\qquad\mbox{on  }\Gamma.
\end{equation*}
\hspace*{5cm}\hfill$\square$
\begin{rem}
    Notice that in Theorem \ref{ns:shaped}, the pressure $p'$ is the shape derivative of the pressure $p_t$ which was defined on $\oo_t$.
\end{rem}
The shape derivative $\vy'$ of the solution $\vy$ of the original
Navier--Stokes system (\ref{ns:nonhomo}) is given by
$\tilde\vy'=\vy'$, then we obtain the following corollary by
substituting $\tilde\vy'=\vy'$ and $\tilde\vy=\vy-\vg$ into
(\ref{ns:shape}).
\begin{corollary}
    The shape derivative $\vy'$ of the solution $\vy$ of (\ref{ns:nonhomo}) exists and satisfies the following system
    \begin{equation}\label{ns:shapederivative}
        \left\{
        \begin{array}{ll}
            -\alpha\Delta\vy'+\md\vy'\cdot\vy+\md\vy\cdot\vy'+\n p'=0 &\quad\mbox{in }\oo;\\
            \mdiv\vy'=0&\quad\text{in }\oo;\\
            \vy'=(\md(\vg-\vy)\cdot\vn)\vv_n&\quad\mbox{on }\Gamma.
        \end{array}
        \right.
    \end{equation}
    Moreover, we have $\vy'\in H^1(\mdiv,\oo)$.
\end{corollary}

\subsection{Adjoint state system and gradients of the cost functionals}
This subsection is devoted to the computation of the shape gradients for the cost functionals $J_1(\oo)$ and $J_2(\oo)$ by the adjoint method.

For the cost functional $J_1(\oo)=\int_\oo\frac{1}{2}\seminorm{\vy-\vy_d}^2\dx$, we have
\begin{theorem}\label{thm:a}
    Let $\oo$ be of class $C^2$ and the velocity $\vec V\in \mathrm{E}^2$, the shape gradient $\n J_1$ of the cost functional $J_1(\oo)$ can be expressed as
    \begin{equation}\label{nsa:gradient}
        \n J_1=\left[\frac{1}{2}(\vy-\vy_d)^2+\alpha(\md(\vy-\vg)\cdot\vn)\cdot(\md\vec v\cdot\vn)\right]\vn,
    \end{equation}
    where the adjoint state $\vec v\in H^1_0(\mdiv,\oo)$ satisfies the following linear adjoint system
    \begin{equation}\label{adjoint:a}
\left\{
\begin{array}{lll}
    -\alpha\Delta\vec v-\md\vec v\cdot\vy+{ }^*\md\vy\cdot\vec v+\n q=\vy-\vy_d,&\qquad\mbox{in  }\oo\\
    \mdiv\vec v=0,&\qquad\mbox{in  }\oo\\
    \vec v=0,&\qquad\mbox{on  }\Gamma.
\end{array}
\right.
    \end{equation}
\end{theorem}
\noindent{\bf Proof.}\; Since $J_1(\oo)$ is differentiable with respect to $\vy$, and the state $\vy$ is shape differentiable with respect to $t$, i.e., the shape derivative $\vy'$ exists, we obtain Eulerian derivative of $J_1(\oo)$ with respect to $t$,
\begin{equation}\label{a:b}
    \diff J_1(\oo;\vec V)=\int_\oo (\vy-\vy_d)\cdot\vy'\dx+\int_\Gamma \frac{1}{2}\seminorm{\vy-\vy_d}^2\vec V_n\diff s
\end{equation}
by Hadamard formula (\ref{hadamard}).

By Green formula, we have the following identity
\begin{multline}\label{a:a}
    \int_\oo [(-\alpha\Delta\vy'+\md\vy'\cdot\vy+\md\vy\cdot\vy'+\n p')\cdot\vw-\mdiv\vy'\pi]\dx\\=\int_\oo[(-\alpha\Delta\vw-\md\vw\cdot\vy+{ }^*\md\vy\cdot\vw+\n\pi)\cdot\vy'-p'\mdiv\vw]\dx\\
+\int_\Gamma (\vy'\cdot\vw)(\vy\cdot\vn)\diff s+\int_\Gamma (\alpha\md\vw\cdot\vn-\pi\vn)\cdot\vy'\diff s+\int_\Gamma (p'\vn-\alpha\md\vy'\vn)\cdot\vw\diff s.
\end{multline}
Now we define $(\vec v,q)\in H^1_0(\mdiv,\oo)\times L^2(\oo)$ to be the solution of (\ref{adjoint:a}), use (\ref{ns:shapederivative}) and set $(\vw,\pi)=
(\vec v,q)$ in (\ref{a:a}) to obtain
\begin{equation}
\int_\oo (\vy-\vy_d)\cdot\vy'\dx=-\int_\Gamma (\alpha\md\vec v\cdot\vn-q\vn)\cdot\vy'\diff s.
\end{equation}
Since $\vy'=(\md(\vg-\vy)\cdot\vn)\vec V_n$ on the boundary $\Gamma$ and $\mdiv\vy'=0$ in $\oo$, we obtain the Eulerian derivative of $J_1(\oo)$ from (\ref{a:b}),
\begin{equation}
    \diff J_1(\oo;\vec V)=\int_\Gamma \left[\frac{1}{2}\seminorm{\vy-\vy_d}^2+
    \alpha\left(\md(\vy-\vg)\cdot\vn\right)\cdot(\md\vec v\cdot\vn)\right]\vec V_n\diff s.
\end{equation}
Since the mapping $\vec V\mapsto \diff J_1(\oo;\vec V)$ is linear and continuous, we get the expression (\ref{nsa:gradient}) for the shape gradient $\n J_1$ by (\ref{pri:shaped}).\hfill $\square$

For another typical cost functional $J_2(\oo)=\frac{\alpha}{2}\int_\oo\seminorm{\mcurl\vy}^2\dx$, we have the following theorem.
\begin{theorem}
    Let $\oo$ be of class $C^2$ and the velocity $\vec V\in \mathrm{E}^2,$ the cost functional $J_2(\oo)$ possesses the shape gradient $\n J_2$ which can be expressed as
    \begin{equation}\label{nsb:gradient}
        \n J_2=\alpha\left[\frac{1}{2}\seminorm{\mathrm{curl}\,\vy}^2+(\md(\vy-\vg)\cdot\vn)\cdot(\md\vec v\cdot\vn-\mathrm{curl}\,\vy\wedge\vn)\right]\vn,
    \end{equation}
    where the adjoint state $\vec v\in H^1_0(\mdiv,\oo)$ satisfies the following linear adjoint system
    \begin{equation}\label{adjoint:b}
\left\{
\begin{array}{lll}
    -\alpha\Delta\vec v-\md\vec v\cdot\vy+{ }^*\md\vy\cdot\vec v+\n q=-\alpha\Delta\vy,&\qquad\mbox{in  }\oo\\
    \mdiv\vec v=0,&\qquad\mbox{in  }\oo\\
    \vec v=0,&\qquad\mbox{on  }\Gamma.
\end{array}
\right.
    \end{equation}
\end{theorem}
\noindent{\bf Proof.}\; The proof is similar to that of Theorem \ref{thm:a}. Using Hadamard formula (\ref{hadamard}) for the cost functional $J_2$, we obtain the Eulerian derivative
\begin{equation}\label{b:b}
    \diff J_2(\oo;\vec V)=\alpha\int_\oo\mcurl\vy\cdot\mcurl\vy'\dx+\int_\Gamma\frac{\alpha}{2}\seminorm{\mcurl\vy}^2\vec V_n\diff s.
\end{equation}
Then, we define $(\vec v,q)\in H^1_0(\mdiv,\oo)\times L^2(\oo)$ to be the solution of (\ref{adjoint:b}), use (\ref{ns:shapederivative}) and set $(\vw,\pi)=
(\vec v,q)$ in (\ref{a:a}) to obtain
\begin{equation}\label{b:c}
\alpha\int_\oo\Delta\vy\cdot\vy'\dx=\int_\Gamma \alpha(\md\vec v\cdot\vn)\cdot\vy'\diff s.
\end{equation}
Applying the following vectorial Green formula
\begin{multline*}
\int_\oo (\vphi\cdot\Delta\vpsi+\mcurl\vphi\cdot\mcurl\vpsi+\mdiv\vphi\mdiv\vpsi)\dx\\=\int_\Gamma (\vphi\cdot(\mcurl\vpsi\wedge\vn)+\vphi\cdot\vn\mdiv\vpsi)\diff s
\end{multline*}
for the vector functions $\vy\in H^1(\mdiv,\oo)$ and $\vy'\in H^1(\mdiv,\oo)$, we obtain
\begin{equation}\label{b:d}
\alpha\int_\oo\mcurl\vy\cdot\mcurl\vy'\dx+\alpha\int_\oo\Delta\vy\cdot\vy'\dx=\alpha\int_\Gamma (\mcurl\vy\wedge\vn)\cdot\vy'\diff s
\end{equation}
Combining (\ref{b:b}), (\ref{b:c}) with (\ref{b:d}), we obtain the Eulerian derivative
\begin{equation}
\diff J_2(\oo;\vec V)=\int_\Gamma\alpha\left[\frac{1}{2}\seminorm{\mathrm{curl}\,\vy}^2+(\md(\vy-\vg)\cdot\vn)\cdot(\md\vec v\cdot\vn-\mathrm{curl}\,\vy\wedge\vn)\right]\vec V_n\diff s.
\end{equation}
Finally we arrive at the expression (\ref{nsb:gradient}) for the shape gradient $\n J_2$.\hfill $\square$

\section{Function space parametrization and function space embedding}\label{sec4}
In this section, we restrict our study to the minimization problem
(\ref{ns:cost}), and problem (\ref{ns:cost2}) follows similarly. In
section \ref{sec:stated}, we have used the local differentiability of the
state with respect to the shape of the fluid domain and the
associated adjoint system to derive the shape gradient of the given
cost functional. However, we do not need to analyze the
differentiability of the state in many cases. In this section we
derive the structure of the shape gradient for the cost functional
$J_1(\oo)=\frac{1}{2}\int_\oo\seminorm{\vy-\vyd}^2\dx$ by function
space parametrization and function space embedding techniques in
order to bypass the study of the state derivative.
\subsection{A saddle point formulation}
In this subsection, we shall describe how to build an appropriate Lagrange functional that takes account into the divergence condition and the nonhomogeneous Dirichlet boundary condition.

We set $\vf\in [H^1(\rn)]^N$ and $\vg\in H^{5/2}(\mdiv,\rn)$, then introduce a Lagrange multiplier $\vec\mu$ and a functional
\begin{equation}
  L(\oo,\vy,p,\vvv,q,\vec\mu)=\int_\oo [(\alpha\Delta\vy-\md\vy\cdot\vy-\n p+\vf)\cdot\vvv+\mdiv\vy q]\dx+\int_\G(\vy-\vg)\cdot\vec\mu\ds
\end{equation}for
$(\vy,p)\in Y(\oo)\times Q(\oo)$, $(\vvv,q)\in P(\oo)\times Q(\oo)$,
 and $\vec\mu\in
H^{-1/2}(\G)^N$ with
\begin{equation*}
Y(\oo):= H^2(\oo)^N; \qquad P(\oo):=
H^2(\oo)^N\cap H^1_0(\oo)^N;\quad Q(\oo):= H^1(\oo).
\end{equation*}
Now we're interested in the following saddle point problem
\begin{equation*}
  \inf_{(\vy,p)\in Y(\oo)\times Q(\oo)}\quad\sup_{(\vvv,q,\vec\mu)\in P(\oo)\times Q(\oo)\times H^{-1/2}(\G)^N}\;L(\oo,\vy,p,\vvv,q,\vec\mu)
\end{equation*}
The solution is characterized by the following systems:
\begin{itemize}
    \item [(i)] The state $(\vy,p)$ is the solution of the problem
\begin{equation}
  \left\{%
  \begin{array}{ll}
-\alpha\Delta\vy+\md\vy\cdot\vy+\n p=\vf&\quad\mbox{in}\;\oo\\
\mdiv \vy=0&\quad\mbox{in}\;\oo\\
\vy=\vg&\quad\mbox{on}\;\G
  \end{array}%
  \right.
\end{equation}
\item [(ii)]The adjoint state $(\vvv,q)$ is the solution of the problem
\begin{equation}\label{prob:bvp1}
\left\{%
  \begin{array}{ll}
      -\alpha\Delta \vvv-\md\vvv\cdot\vy+{ }^*\md\vy\cdot\vvv+\n q=0\qquad&\mbox{in}\;\oo\\
\mdiv \vvv=0&\mbox{in}\;\oo\\
    \vvv=0\qquad&\mbox{on}\;\G;
  \end{array}%
  \right.
\end{equation}

 \item [(iii)] The multiplier: $\vec\mu=\alpha\md\vvv\,\vec{n}-q\,\vn,\;\mbox{on}\;\G$.
\end{itemize}
Hence we obtain the following new functional,
\begin{equation*}
  L(\oo,\vy,p,\vvv,q)=\int_\oo [(\alpha\Delta\vy-\md\vy\cdot\vy-\n p+\vf)\cdot\vvv+\mdiv\vy q]\dx+\int_\G(\vy-\vg)\cdot(\alpha\md\vvv\cdot\vec{n}-q\vn)\ds.
\end{equation*}
To get rid of the boundary integral, the following identities are
derived by Green formula,
\begin{eqnarray*}
\int_\G(\vy-\vg)\cdot(\md\vvv\cdot\vn)\ds&=&\int_\oo[(\vy-\vg)\cdot\Delta\vvv+\md(\vy-\vg):\md\vvv]\dx;\\
  \int_\G(\vy-\vg)\cdot(q\,\vec{n})\ds&=&\int_\oo
  [\mdiv(\vy-\vg)q+(\vy-\vg)\cdot\n q]\dx.
\end{eqnarray*}
Thus we introduce the new Lagrangian associated with
(\ref{ns:nonhomo}) and the cost functional
$J_1(\oo)=\frac{1}{2}\int_\oo\seminorm{\vy-\vyd}^2\dx$:
\begin{multline*}
  G(\oo,\vy,p,\vvv,q)=\frac{1}{2}\int_\oo\seminorm{\vy-\vyd}^2\dx+\int_\oo [(\alpha\Delta\vy-\md\vy\cdot\vy-\n p+\vf)\cdot\vvv+\mdiv\vy q]\dx\\+
  \alpha\int_\oo[(\vy-\vg)\cdot\Delta\vvv+\md(\vy-\vg):\md\vvv]\dx-\int_\oo
  [\mdiv(\vy-\vg)q+(\vy-\vg)\cdot\n q]\dx.
\end{multline*}
Now the minimization problem (\ref{ns:cost}) can be expressed as the following form
\begin{equation*}
    \min_{\oo\in\mathcal{O}} \inf_{(\vy,p)\in Y(\oo)\times P(\oo)}\sup_{(\vvv,q)\in P(\oo)\times Q(\oo)}G(\oo,\vy,p,\vvv,q).
\end{equation*}
\indent We can use the minimax framework to avoid the study of the
state derivative with respect to the shape of the domain. The
Karusch-Kuhn-Tucker (KKT) conditions will furnish the shape gradient
of the cost functional $J_1(\oo)$ by using the adjoint system. To
begin with, we derive the formulation of the adjoint system which
was satisfied by $(\vvv,q)$.

\indent For $(p,\vvv,q)\in Q(\oo)\times P(\oo)\times Q(\oo)$,
$G(\oo,\vy,p,\vvv,q)$ is differentiable with respect to $\vy\in
Y(\oo)$ and we get
\begin{multline*}
        {\p_{\vy}} G(\oo,\vy,p,\vvv,q)\cdot\delta\vy=\int_\oo (\vy-\vyd)\cdot\delta\vy\dx+\int_\oo [\alpha\Delta(\delta\vy)-\md(\delta\vy)\cdot\vy-\md\vy\cdot\delta\vy]\cdot\vvv\dx\\+\alpha\int_\oo[\delta\vy\cdot\Delta\vvv+\md(\delta\vy):\md\vvv]\dx-\int_\oo \delta\vy\cdot\n q\dx,\quad\forall\delta\vy\in Y(\oo).
\end{multline*}
Integrating by parts, we obtain
\begin{equation}
        \label{kkt:y}
        {\p_{\vy}} G(\oo,\vy,p,\vvv,q)\cdot\delta\vy=\int_\oo(\alpha\Delta\vvv+\md\vvv\cdot\vy-{ }^*\md\vy\cdot\vvv-\n q+\vy-\vyd)\cdot\delta\vy\dx.
\end{equation}
Similarly for $(\vy,\vvv,q)\in Y(\oo)\times P(\oo)\times Q(\oo)$, $G(\oo,\vy,p,\vvv,q)$ is differentiable with respect to $p\in Q(\oo)$, and we have
\begin{equation}
        \label{kkt:p}
        {\p_p}G(\oo,\vy,p,\vvv,q)\cdot\delta p=\int_\oo-\n(\delta p)\cdot\vvv\dx=\int_\oo\delta p\,\mdiv\vvv\dx,\quad\forall\delta p\in Q(\oo).
\end{equation}
Hence, (\ref{kkt:y}) and (\ref{kkt:p}) lead to the following linear adjoint system
\begin{equation}\label{fse:adjoint}
\left\{%
  \begin{array}{ll}
      -\alpha\Delta \vvv-\md\vvv\cdot\vy+{ }^*\md\vy\cdot\vvv+\n q=\vy-\vyd\qquad&\mbox{in}\;\oo\\
\mdiv \vvv=0&\mbox{in}\;\oo\\
    \vvv=0\qquad&\mbox{on}\;\G;
  \end{array}%
  \right.
\end{equation}

Given a velocity field $\vec V\in \mathrm{E}^2$ and transformed
domain $\oo_t:=T_t(\oo)$, our main task of this section is to get the limit
\begin{equation}
  \lim_{t\searrow 0}\frac{j(t)-j(0)}{t}
\end{equation}
with
\begin{equation}\label{lf:jt}
j(t)= \inf_{(\vy_t,p_t)\in Y(\oo_t)\times
Q(\oo_t)}\quad\sup_{(\vvv_t,q_t)\in P(\oo_t)\times
Q(\oo_t)}G(\oo_t,\vy_t,p_t,\vvv_t,q_t),
\end{equation}
 where $(\vy_t,p_t)$ satisfies
\begin{equation}
  \left\{%
  \begin{array}{ll}
-\alpha\Delta\vy_t+\md\vy_t\cdot\vy_t+\n p_t=\vf&\quad\mbox{in}\;\oo_t\\
\mdiv \vy_t=0&\quad\mbox{in}\;\oo_t\\
\vyt=\vg&\quad\mbox{on}\;\G_t
  \end{array}%
  \right.
\end{equation}
and $(\vec v_t,q_t)$ satisfies
\begin{equation}
\left\{%
  \begin{array}{ll}
    -\alpha\Delta \vvv_t-\md\vvv_t\cdot\vy_t+{ }^*\md\vy_t\cdot\vvv_t+\n q_t=\vy_t-\vyd\qquad&\mbox{in}\;\oo_t\\
\mdiv \vvv_t=0&\mbox{in}\;\oo_t\\
    \vvv_t=0\qquad&\mbox{on}\;\G_t;
  \end{array}%
  \right.
\end{equation}
 Unfortunately, the Sobolev space $Y(\oo_t)$, $Q(\oo_t)$, and
$P(\oo_t)$ depend on the parameter $t$, so we need a theorem to
differentiate a saddle point with respect to the parameter $t$, and
there are two techniques to get rid of it:
\begin{itemize}
    \item   {Function space parametrization }technique;
    \item { Function space embedding }technique.
\end{itemize}
\subsection{Function space parametrization}\label{fsp}
This subsection is devoted to the {function space
parametrization}, which consists in transporting the different
quantities (such as, a cost functional) defined on the variable domain
$\oo_t$ back into the reference domain $\oo$ which does not depend
on the perturbation parameter $t$. Thus we can use differential
calculus since the functionals involved are defined in a fixed
domain $\oo$ with respect to the parameter $t$.

We parameterize the functions in $H^m(\oo_t)^d$ by elements of
$H^m(\oo)^d$ through the transformation:
\begin{equation*}
  \vphi\mapsto \vphi\circ T_t^{-1}:\quad H^m(\oo)^d\rightarrow
  H^m(\oo_t)^d,\qquad \mbox{integer}\;m\geq 0.
\end{equation*}
where "$\circ$" denotes the composition of the two maps and $d$ is
the dimension of the function $\vphi$.

 Note that
since $T_t$ and $T_t^{-1}$ are diffeomorphisms, it transforms the
reference domain $\oo$ (respectively, the boundary $\G$) into the
new domain $\oo_t$ (respectively, the boundary $\G_t$ of $\oo_t$).
This parametrization can not change the value of the saddle point.
We can rewrite (\ref{lf:jt}) as
\begin{equation}\label{fsp:newsaddlep}
j(t)= \inf_{(\vy,p)\in Y(\oo)\times
Q(\oo)}\quad\sup_{(\vvv,q)\in P(\oo)\times Q(\oo)}G(\oo_t,\vy\circ
T_t^{-1},p\circ T_t^{-1},\vvv\circ T_t^{-1},q\circ T_t^{-1}).
\end{equation}
It amounts to introducing the new Lagrangian for $(\vy,p,\vvv,q)\in
Y(\oo)\times Q(\oo)\times P(\oo)\times Q(\oo)$:
\begin{equation*}
  \tilde G(t,\vy,p,\vvv,q)\defmath G(\oo_t,\vy\circ
T_t^{-1},p\circ T_t^{-1},\vvv\circ T_t^{-1},q\circ T_t^{-1}).
\end{equation*}
The expression for $\tilde G(t,\vy,p,\vvv,q)$ is given by
\begin{equation}\label{gtp}
  \tilde G(t,\vy,p,\vvv,q):=I_1(t)+I_2(t)+I_3(t)+I_4(t),
\end{equation} where
\begin{eqnarray*}
 I_1(t)&:=&\frac{1}{2}\int_{\oo_t}\seminorm{\vy\circ
 T_t^{-1}-\vyd}^2\dx;\\
 I_2(t)&:=&\int_{\oo_t} [(\alpha\Delta(\vy\circ
 T_t^{-1})-\md(\vy\circ T_t^{-1})\cdot(\vy\circ T_t^{-1})-\n (p\circ T_t^{-1})+\vf]\cdot(\vvv\circ T_t^{-1})\\
 &&\quad\quad\qquad+\mdiv(\vy\circ
T_t^{-1})(q\circ T_t^{-1})\dx;\\
I_3(t)&:=&\alpha\int_{\oo_t}[(\vy\circ
T_t^{-1}-\vg\cdot\Delta(\vvv\circ T_t^{-1})+\md(\vy\circ
T_t^{-1}-\vg):\md(\vvv\circ
T_t^{-1})]\dx;\\
I_4(t)&:=&-\int_{\oo_t}[\mdiv(\vy\circ T_t^{-1}-\vg)(q\circ
T_t^{-1})+(\vy\circ T_t^{-1}-\vg)\cdot\n(q\circ
T_t^{-1})]\dx,
\end{eqnarray*}
Now we introduce the theorem concerning on the differentiability of
a saddle point (or a minimax). To begin with, some notations are
given as follows.

 Define a functional
$$\mg : [0,\tau]\times X\times Y\rightarrow\mathbb{R}$$
with $\tau>0$, and $X,Y$ are the two topological spaces.

 For any
$t\in [0,\tau]$, define
$$g(t)=\inf_{x\in X}\sup_{y\in Y}\mg(t,x,y)$$
and the sets
\begin{eqnarray*}
  &X(t)=\{x^t\in X:g(t)=\sup_{y\in Y}\mg(t,x^t,y)\}\\
  &Y(t,x)=\{y^t\in Y:\mg(t,x,y^t)=\sup_{y\in Y}\mg(t,x,y)\}
\end{eqnarray*}
Similarly, we can define dual functionals
$$h(t)=\sup_{y\in Y}\inf_{x\in X}\mg(t,x,y)$$
and the corresponding sets
\begin{eqnarray*}
 & Y(t)=\{y^t\in Y:h(t)=\inf_{x\in X}\mg(t,x,y^t)\}\\
  &X(t,y)=\{x^t\in X:\mg(t,x^t,y)=\inf_{x\in X}\mg(t,x,y)\}
\end{eqnarray*}
Furthermore, we introduce the set of saddle points
$$S(t)=\{(x,y)\in X\times Y: g(t)=\mg(t,x,y)=h(t)\}$$
Now we can introduce the following theorem (see \cite{correa} or
page 427 of \cite{delfour}):
\begin{theorem}\label{fsp:correa}
 Assume that the following hypothesis hold:
 \begin{itemize}
    \item [(H1)]$S(t)\neq\emptyset,\;t\in [0,\tau];$
    \item [(H2)]The partial derivative $\p_t\mg(t,x,y)$ exists in
    $[0,\tau]$ for all $$(x,y)\in \left[\underset{{t\in [0,\tau]}}{\bigcup}X(t)\times Y(0)\right]\bigcup\left[X(0)\times\underset{{t\in [0,\tau]}}{\bigcup}Y(t)\right];$$
    \item [(H3)]There exists a topology $\mt_X$ on $X$ such that for
    any sequence $\{t_n:t_n\in [0,\tau]\}$ with
    $\lim\limits_{n\nearrow\infty}t_n=0$, there exists $x^0\in X(0)$ and a subsequence
    $\{t_{n_k}\}$, and for each $k\geq 1,$ there exists $x_{n_k}\in
    X(t_{n_k})$ such that
    \begin{enumerate}
        \item [(i)]$\lim\limits_{n\nearrow\infty}x_{n_k}=x^0$ in the
        $\mt_X$ topology,
        \item [(ii)]$$\liminf\limits_{t\searrow 0\atop k\nearrow\infty}\p_t\mg(t,x_{n_k},y)\geq\p_t\mg(0,x^0,y),\quad \forall y\in Y(0);$$
    \end{enumerate}
    \item [(H4)]There exists a topology $\mt_Y$ on $Y$ such that for
    any sequence $\{t_n:t_n\in [0,\tau]\}$ with
    $\lim\limits_{n\nearrow\infty}t_n=0$, there exists $y^0\in Y(0)$ and a subsequence
    $\{t_{n_k}\}$, and for each $k\geq 1,$ there exists $y_{n_k}\in
    Y(t_{n_k})$ such that
    \begin{enumerate}
        \item [(i)]$\lim\limits_{n\nearrow\infty}y_{n_k}=y^0$ in the
        $\mt_Y$ topology,
        \item [(ii)]$$\limsup\limits_{t\searrow 0\atop k\nearrow\infty}
        \p_t\mg(t,x,y_{n_k})\leq\p_t\mg(0,x,y^0),\quad \forall x\in X(0).$$
    \end{enumerate}
 \end{itemize}
 Then there exists $(x^0,y^0)\in X(0)\times Y(0)$ such that
 \begin{multline}
   \diff g(0)=\lim_{t\searrow 0}\frac{g(t)-g(0)}{t}\\=\inf_{x\in
   X(0)}\sup_{y\in Y(0)}\p_t \mg(0,x,y)=\p_t\mg(0,x^0,y^0)=\sup_{y\in Y(0)}\inf_{x\in
   X(0)}\p_t \mg(0,x,y)
 \end{multline}
 This means that $(x^0,y^0)\in X(0)\times Y(0)$ is a saddle point of
 $\p_t\mg(0,x,y)$.
\end{theorem}
Following Theorem \ref{fsp:correa}, we need to differentiate the perturbed Lagrange functional $\tilde G(t,\vy,p,\vvv,q)$. Since $(\vy,p,\vvv,q)\in H^3(\oo)^N\times H^2(\oo)\times H^3(\oo)^N\times
H^2(\oo)$ provided that $\G$ is at less $C^3$ (see \cite{temam01}),
 we can use Hadamard formula (\ref{hadamard}) to
 differentiate $\tilde G(t,\vy,p,\vvv,q)$ with respect to the parameter $t>0$,
$$\p_t\tilde G(t,\vy,p,\vvv,q)=I'_1(0)+I'_2(0)+I'_3(0)+I'_4(0),$$
where
\begin{equation}\label{i1} I'_1(0):=\frac{\p}{\p
t}\left\{I_1(t)\right\}\Big{|}_{t=0}=\int_\oo
(\vy-\vyd)\cdot(-\md\vy\vec
  V)\dx+\frac{1}{2}\int_\G\seminorm{\vy-\vyd}^2\vec
  V_n\ds;
\end{equation}
\begin{multline}\label{i2}
I'_2(0):=\frac{\p}{\p
t}\left\{I_2(t)\right\}\Big{|}_{t=0}=\int_\oo(\alpha\Delta\vy-\md\vy\cdot\vy-\n
p+\vf)\cdot(-\md\vvv\vec
V)\dx\\
+\int_{\oo}
[\alpha\Delta(-\md\vy\vec V)-\md(-\md\vy\vv)\cdot\vy-\md\vy\cdot(-\md\vy\vv)+\n(\n p\cdot\vec
V)]\cdot\vvv\dx\\+\int_\oo [\mdiv(-\md\vy\vec V)q+\mdiv\vy(-\n q\cdot\vec V)]\dx\\
+\int_\G[(\alpha\Delta\vy-\md\vy\cdot\vy-\n
p+\vf)\cdot\vvv+\mdiv\vy\,q]\vec V_n\ds;
\end{multline}
\begin{multline}\label{i3}
I'_3(0):=\frac{\p}{\p
t}\left\{I_3(t)\right\}\Big{|}_{t=0}=\alpha\int_\oo
\left\{(-\md\vy\vec
V)\cdot\Delta\vvv+(\vy-\vg)\cdot\Delta(-\md\vvv\vec
  V)\right.\\\left.{\qquad}+\md(-\md\vy\vec V):\md\vvv+\md(\vy-\vg):\md(-\md\vvv\vec
  V)\right\}\dx\\+\alpha\int_\G [(\vy-\vg)\cdot\Delta\vvv\rangle+\md(\vy-\vg):\md\vvv]\vec
  V_n\ds;
\end{multline}
\begin{multline}\label{i4}
I'_4(0):=\frac{\p}{\p t}\left\{I_4(t)\right\}\Big{|}_{t=0}=-\int_\oo
[\mdiv(-\md\vy\vec V)q+\mdiv(\vy-\vg)(-\n q\cdot\vec
V)\\
{\hspace*{1cm}}+(-\md\vy\vec V)\cdot\n q-(\vy-\vg)\cdot\n(\n
q\cdot\vec V)]\dx-\int_\G [\mdiv(\vy-\vg)q+(\vy-\vg)\cdot\n
q]\vec V_n\ds.
\end{multline}
Since $(\vy,p)$ satisfies (\ref{ns:nonhomo}), $\mdiv \vec v=0$ and $\vec v|_\G=0$, also by Green formula we can simplify (\ref{i2}) to
\begin{multline}\label{ii2}
I'_2(0)=\int_{\oo}
[\alpha\Delta(-\md\vy\vec V)-\md(-\md\vy\vv)\cdot\vy-\md\vy\cdot(-\md\vy\vv)]\cdot\vvv\dx\\
+\int_\oo\mdiv(-\md\vy\vec V)q\dx.
\end{multline}
By Green formula, $\vy|_\G=\vg$ and $\vec v|_\G=0$, we can simplify (\ref{i3}) to
\begin{multline}\label{ii3}
I'_3(0)=\alpha\int_\oo (-\md\vy\vec
V)\cdot\Delta\vvv\dx-\alpha\int_\oo\Delta(-\md\vy\vec
V)\cdot\vvv\dx\\+\alpha\int_\G [\md(\vy-\vg):\md\vvv]\vec
  V_n\ds;
\end{multline}
By Green formula, $\vy|_\G=\vg$ and $\mdiv\vy=\mdiv\vg=0$, (\ref{i4}) can be simplified to
\begin{equation}\label{ii4}
I'_4(0)=-\int_\oo [(-\md\vy\vec V)\cdot\n q+\mdiv(-\md\vy\vec
V)q]\dx
\end{equation}
Adding (\ref{i1}), (\ref{ii2}), (\ref{ii3}) and (\ref{ii4}) together,
\begin{multline}
    \label{ii5}
    \sum^4_{i=1}I'_i(0)\\
=-\int_\oo[\alpha\Delta\vvv\cdot(\md\vy\vv)-\md(\md\vy\vv)\cdot\vy-\md\vy\cdot(\md\vy\vv)-(\md\vy\vv)\cdot\n q+(\vy-\vyd)\cdot(\md\vy\vv)]\dx\\
+\int_\G\left\{\frac{1}{2}\seminorm{\vy-\vyd}^2+\alpha\md(\vy-\vg):\md\vvv
\right\}\vec V\cdot\vn\ds.
\end{multline}
Since $(\vvv,q)$ are characterized by (\ref{fse:adjoint}), we
multiply the first equation of (\ref{fse:adjoint}) by $(\md\vy\vv)$
and integrate over $\oo$, then the distributional integral in
(\ref{ii5}) vanishes, finally we obtain the boundary expression for
the Eulerian derivative
\begin{equation}\label{end}
 \diff J_1(\oo;\vec V)=\int_\G\left\{\frac{1}{2}\seminorm{\vy-\vyd}^2+\alpha\md(\vy-\vg):\md\vvv
\right\}\vec V\cdot\vn\ds.
\end{equation}
Since $\vy|_\G=\vg$ and $\vec v|_\G=0$, we
have $\md(\vy-\vg)|_\G=\md(\vy-\vg)\cdot\vn{ }^*\vn$ and $\md\vvv|_\G=\md\vvv\cdot\vn{ }^*\vn$,
 thus we obtain an expression for the shape gradient
\begin{equation}
        \label{fsp:gradient}
        \n
J_1=\left\{\frac{1}{2}\seminorm{\vy-\vyd}^2+\alpha[\md(\vy-\vg)\cdot\vn]\cdot(\md\vvv\cdot\vn)\right\}\vn,
\end{equation}
which is the same as (\ref{nsa:gradient}) in the previous section.
\subsection{Function space embedding}\label{fse}
In the previous subsection, we have used the technique of function
space parametrization in order to get the derivative of $j(t)$,
i.e.,
\begin{equation}
j(t)=\inf_{(\vy_t,p_t)\in Y(\oo_t)\times Q(\oo_t)}\sup_{(\vvv_t,q_t)\in
  P(\oo_t)\times Q(\oo_t)}G(\oo_t,\vy_t,p_t,\vvv_t,q_t).
\end{equation}
with respect to the parameter
$t>0.$ This subsection is devoted to a different method based on
function space embedding technique. It means that the state and
adjoint state are defined on a large enough domain $D$ (called a
\textit{hold-all} \cite{delfour}) which contains all the
transformations $\{\oo_t: 0\leq t\leq\varepsilon\}$ of the reference domain
$\oo$ for some small $\varepsilon>0.$

For convenience, let $D=\rn$. Use the function space embedding technique,
\begin{equation}\label{fse:fun}
  j(t)=\inf_{(\vec\my,\mpp)\in Y(\rn)\times Q(\rn)}\sup_{(\vec\mv,\mq)\in
  P(\rn)\times Q(\rn)}G(\oo_t,\vec\my,\mpp,\vec\mv,\mq).
\end{equation}
where the new Lagrangian
\begin{multline}\label{fse:lag}
G(\oo_t,\vec{\mathcal Y},{\mathcal P},\vec\mv,\mathcal{Q})=\frac{1}{2}\int_{\oo_t}\seminorm{\vec\my-\vyd}^2\dx
                    +\int_{\oo_t}[(\alpha\Delta\vec\my-\md\vec\my\cdot\vec\my-\n \mpp+\vf)\cdot\vec\mv+\mdiv\vec\my\, \mq]\dx\\
  +
  \alpha\int_{\oo_t}[(\vec\my-\vg)\cdot\Delta\vec\mv+\md(\vec\my-\vg):\md\vec\mv]\dx-\int_{\oo_t}
  [\mdiv(\vec\my-\vg)\mq+(\vec\my-\vg)\cdot\n \mq]\dx.
\end{multline}
Since $\vf,\vyd\in H^1(\rn)^N,$ $\vg\in H^{5/2}(\rn)^N$, and $\oo_t$
is of class $C^3$, the solution $(\vy_t,p_t,\vvv_t,q_t)$ of (\dots)
belongs to $H^3(\oo_t)^N\times (H^2(\oo_t)\cap L^2_0(\oo_t))\times
(H^3(\oo_t)^N\cap H^1_0(\oo_t)^N)\times (H^2(\oo_t)\cap
L^2_0(\oo_t))$ instead of $Y(\oo_t)\times Q(\oo_t)\times
P(\oo_t)\times Q(\oo_t).$ Therefore, the sets
$$X=Y=H^3(\rn)^N\times H^2(\rn),$$
and the saddle points $S(t)=X(t)\times Y(t)$ are given by
\begin{eqnarray}
  X(t)&=&\{(\vec\my,\mpp)\in X: \,\vec\my|_{\oo_t}=\vy_t,\;\mpp|_{\oo_t}=p_t\}\\
  Y(t)&=&\{(\vec\mv,\mq)\in Y: \,\vec\mv|_{\oo_t}=\vvv_t,\;\mq|_{\oo_t}=q_t\}
\end{eqnarray}
Using Theorem \ref{fsp:correa}, we may make the conjecture that we can bypass the inf--sup and state
\begin{equation}\label{dd}
  \diff J(\oo;\vec V)=\inf_{(\vec\my,\mpp)\in X(0)}\sup_{(\vec\mv,\mq)\in Y(0)}\p_t
  G(\oo_t,\vec\my,\mpp,\vec\mv,\mq)|_{t=0}.
\end{equation}
Now we compute the partial derivative of the expression
(\ref{fse:lag}) by Hadamard formula (\ref{hadamard}),
\begin{equation}\label{fse:aaa}
  \p_t
  G(\oo_t,\vec\my,\mpp,\vec\mv,\mq)=\int_{\G_t}[\mathcal{W}_1(\vec\my,\vec\mv)
  +\mathcal{W}_2(\vec\my,\mpp,\vec\mv,\mq)]\vec
  V\cdot\vn_t\ds
\end{equation}
where
\begin{eqnarray*}
  \mathcal{W}_1(\vec\my,\vec\mv)&:=&\frac{1}{2}\seminorm{\vec\my-\vyd}^2+\alpha\md(\vec\my-\vg):\md\vec\mv;\\
 \mathcal{W}_2(\vec\my,\mpp,\vec\mv,\mq)&:=&(\alpha\Delta\vec\my-\md\vec\my\cdot\vec\my-\n \mpp+\vf)\cdot\vec\mv
 +(\vec\my-\vg)\cdot(\alpha\Delta\vec\mv-\n\mq)+\mq\mdiv\vg.
\end{eqnarray*}
and $\vn_t$ denotes the outward unit normal to the boundary $\G_t$.

We note that the expression (\ref{fse:aaa}) is a boundary
integral on $\G_t$ which will not depend on $(\vec\my,\mpp)$ and
$(\vec\mv,\mq)$ outside of $\overline{\oo}_t$, so the inf and the
sup in (\ref{dd}) can be dropped, we then get
\begin{equation*}
  \diff J_1(\oo;\vec V)=\int_\G [\mathcal{W}_1(\vy,\vvv)+\mathcal{W}_2(\vy,p,\vvv,q)]\vec V_n\ds
\end{equation*}
However, $(\vy-\vg)|_\G=\vvv|_\G=0$ and $\mdiv\vg=0$ imply
$\mathcal{W}_2(\vy,p,\vvv,q)=0$ on the boundary $\G$. Finally we
have
\begin{equation*}
  \diff J_1(\oo;\vec V)=\int_\G\left\{\frac{1}{2}\seminorm{\vy-\vyd}^2
  +\alpha\md(\vy-\vg):\md\vvv\right\}\vv_n\ds
\end{equation*}
As in the previous subsection, we also have the shape gradient of the functional $J(\oo)$,
\begin{equation}
    \n
J_1=\left\{\frac{1}{2}\seminorm{\vy-\vyd}^2+\alpha[\md(\vy-\vg)\cdot\vn]\cdot(\md\vvv\cdot\vn)\right\}\vn,
\end{equation}
which is the same as the expressions (\ref{nsa:gradient}) and
(\ref{fsp:gradient}).
\section{Gradient algorithm and numerical simulation}
In this section, we will give a gradient type algorithm and some
numerical examples in two dimensions to prove that our previous
methods could be very useful and efficient for the numerical
implementation of the shape optimization problems for Navier--Stokes
flow.
 \subsection{A gradient type algorithm}
 In this subsection, we will describe a gradient
type algorithm for the minimization of a cost function $J(\oo)$. As
we have just seen, the general form of its Eulerian derivative is
\begin{equation*}
  \diff J(\oo;\vec V)=\int_\G \n J\cdot \vec V\ds,
\end{equation*}
where $\n J$ denotes the shape gradient of the cost functional $J$.
Ignoring
regularization, a descent direction is found by defining
\begin{equation}
  \vec V=-h_k\n J
\end{equation}
and then we can update the shape $\oo$ as
\begin{equation}
  \oo_k=(\mathrm{I}+h_k\vec V)\oo
\end{equation}
where $h_k$ is a descent step at $k$-th iteration.

There are also other choices for the definition of the descent
direction. Since the gradient of the functional has necessarily less regularity than the parameter, an iterative scheme like the method of descent deteriortates the regularity of the optimized parameter. We need to project or smooth the variation into $H^1(\oo)^2$. Hence, the method used in this paper is to change the scalar
product with respect to which we compute a descent direction, for
instance, $H^1(\oo)^2$. In this case, the descent direction is the
unique element $\vec d\in H^1(\oo)^2$ such that for every $\vec V\in
H^1(\oo)^2,$
\begin{equation}\label{reg}
  \int_\oo\md\vec d :\md \vec V\dx=\diff J(\oo;\vec V).
\end{equation}
  The computation of $\vec d$ can also be interpreted as a regularization
  of the shape gradient, and the choice of $H^1(\oo)^2$ as space of
  variations is more dictated by technical considerations rather
  than theoretical ones.

The resulting algorithm can be summarized as follows:
\begin{itemize}
    \item [(1)] Choose an initial shape $\oo_0$;
    \item [(2)] Compute the state system and adjoint state system, then
    we can evaluate the descent direction $\vec d_k$ by using (\ref{reg})
    with $\oo=\oo_k;$
    \item[(3)] Set $\oo_{k+1}=(\mathrm{Id}-h_k\vec d_k) \,\oo_k,$ where $h_k$
    is a small positive real number.
\end{itemize}

The choice of the descent step $h_k$ is not an easy task. Too big, the algorithm is unstable; too small, the rate of convergence is insignificant. In order to refresh $h_k$, we compare $h_k$ with $h_{k-1}$. If $(\vec d_k,\vec d_{k-1})_{H^1}$ is negative, we should reduce the step; on the other hand, if $\vec d_k$ and $\vec d_{k-1}$ are very close, we increase the step. In addition, if reversed triangles are appeared when moving the mesh, we also need to reduce the step.

In our algorithm, we do not choose any stopping criterion. A classical stopping criterion is to find that whether the shape gradients $\n J$ in some suitable norm is small enough. However, since we use the continuous shape gradients, it's hopeless for us to expect very small gradient norm because of numerical discretization errors. Instead, we fix the number of iterations. If it is too small, we can restart it with the previous final shape as the initial shape.
\subsection{Numerical examples}
 To illustrate the theory, we
want to solve the following minimization problem
\begin{equation}\label{exam:fun}
    \min_{\oo}\frac{1}{2}\int_{\oo}\seminorm{\vy-\vyd}^2\dx
\end{equation}
subject to
\begin{equation}
  \left\{%
  \begin{array}{ll}
-\alpha\Delta\vy+\md\vy\cdot\vy+\n p=\vf&\quad\mbox{in}\;\oo\\
\mdiv \vy=0&\quad\mbox{in}\;\oo\\
\vy=0&\quad\mbox{on}\;\G:=\G_{\mathrm{out}}\cup\G_{\mathrm{in}};\\
  \end{array}%
  \right.
\end{equation}
The domain $\oo$ is an annuli, and its boundary has two parts: the
outer boundary $\G_{\mathrm{out}}$ is an unit circle which is fixed;
the inner boundary $\G_{\mathrm{in}}$ which is to be optimized. We
choose the body force $\vf=(f_1,f_2)$ as follows:
\begin{multline*}
 f_1=-x^3+\frac{\alpha y (15 x^2+15 y^2-1)}{5(x^2+y^2)^{3/2}}\\-x\left[y^2+\frac{1}{775}
  \left(1426+\frac{31}{x^2+y^2}+\frac{753}{\sqrt{x^2+y^2}}-1860\sqrt{x^2+y^2}\right)\right];
\end{multline*}
\begin{multline*}
 f_2=-y^3-\frac{\alpha y (15 x^2+15 y^2-1)}{5(x^2+y^2)^{3/2}}\\-y\left[x^2+\frac{1}{775}
  \left(1426+\frac{31}{x^2+y^2}+\frac{753}{\sqrt{x^2+y^2}}-1860\sqrt{x^2+y^2}\right)\right].
\end{multline*}
The target velocity $\vyd$ is determined by the data $\vf$ and the target shape of the domain $\oo$.

In this model problem, we have the following Eulerian derivative:
\begin{equation*}
\diff J(\oo;\vv)=\int_{\G_{\mathrm{in}}}\left\{\frac{1}{2}\seminorm{\vy-\vyd}^2+\alpha\md\vy:\md\vec v\right\}\vec V_n\ds.
\end{equation*}
We will solve this shape problem with two different target shapes:\\[3pt]
\noindent{\bf Case 1:}\; A circle: $\G_{\mathrm{in}}=\{(x,y)|x^2+y^2=0.2^2\}$. \\[4pt]
\noindent{\bf Case 2:}\; An ellipse: $\G_{\mathrm{in}}=\{(x,y)|\frac{x^2}{4}+\frac{4y^2}{9}=\frac{1}{25}\}.$\\
Our numerical solutions are obtained under FreeFem++ \cite{hecht}
and we run the program on a home PC.

 In Case 1, we choose the
initial shape to be elliptic:
$\left\{(x,y)\Big{|}\frac{x^2}{9}+\frac{y^2}{4}=\frac{1}{25}\right\}$,
and the initial mesh is shown in \autoref{fig0}.

In Case 2, we take the initial shape to be a circle whose center is at
origin with
radius 0.6, and the initial mesh is shown in \autoref{fig00}.

\begin{figure}[!htbp]
\renewcommand{\captionlabelfont}{\small}
\setcaptionwidth{1.7in}
\begin{minipage}[b]{0.50\textwidth}
  \centering
  \includegraphics[width=2.1in]{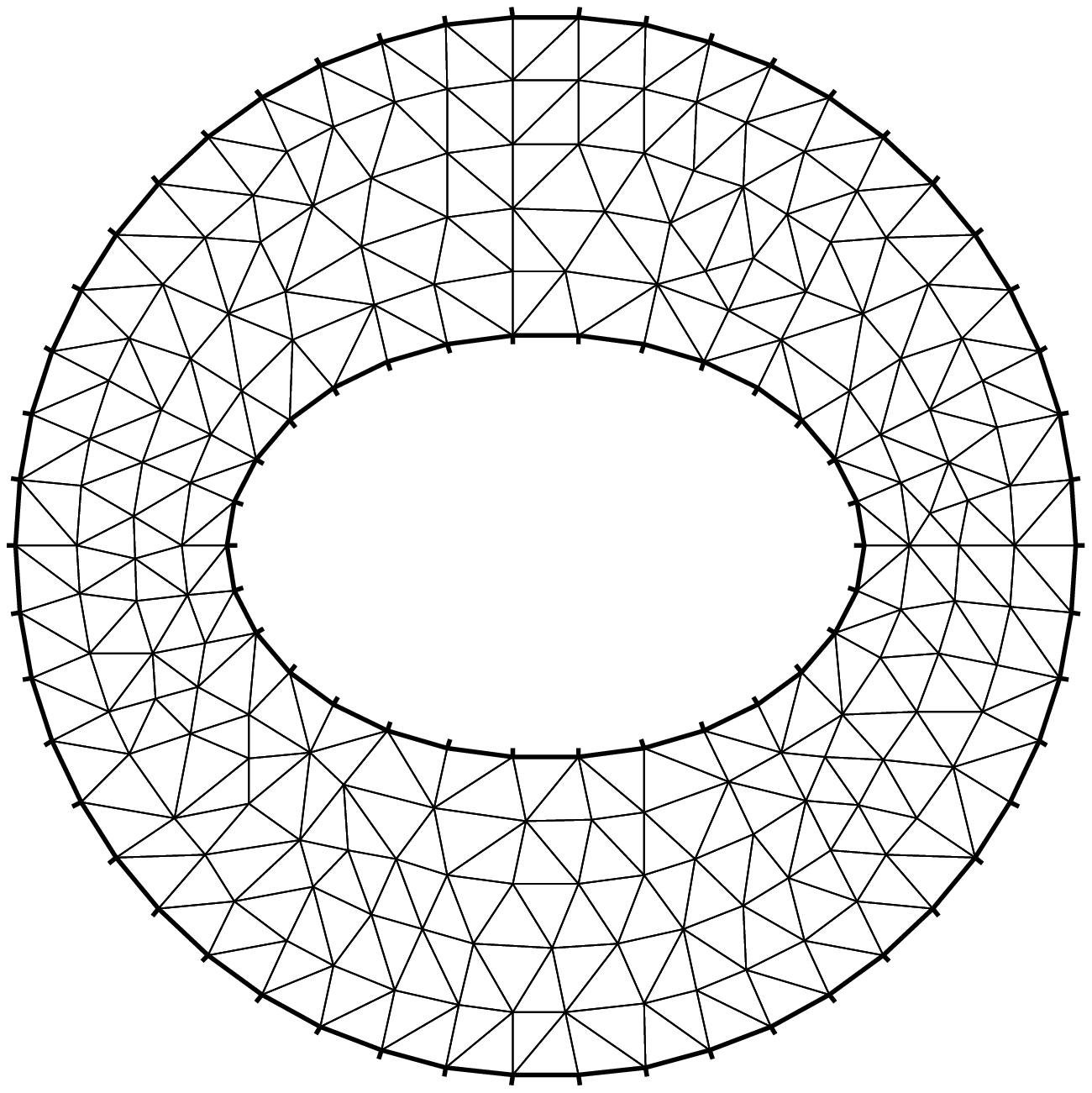}
  \caption{Initial mesh in Case 1 with 226 nodes.\label{fig0}}
  \end{minipage}
\begin{minipage}[b]{0.50\textwidth}
  \centering
  \includegraphics[width=2.1in]{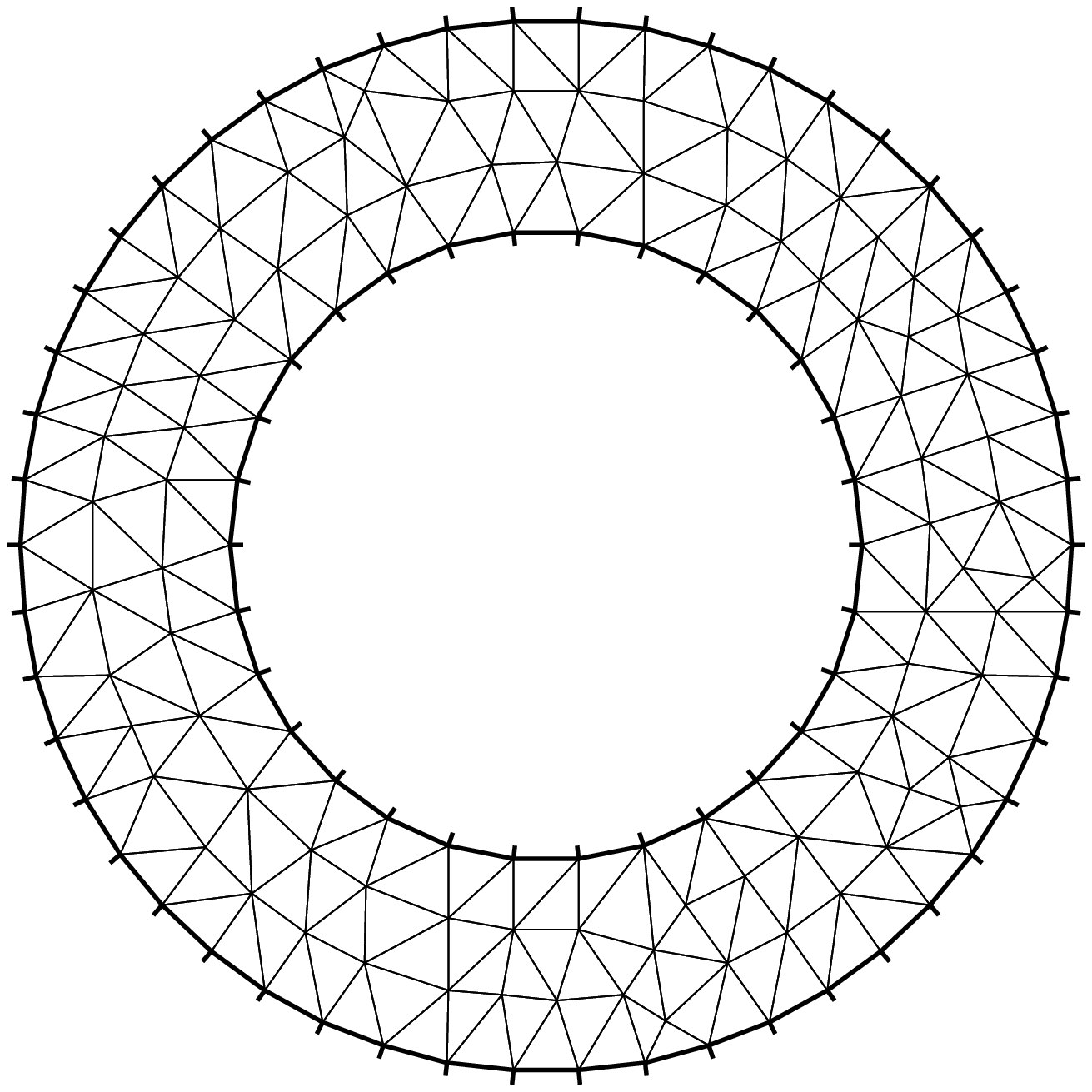}
  \caption{Initial mesh in Case 2 with 161 nodes.\label{fig00}}
  \end{minipage}
\end{figure}

In Case 1,  \autoref{fig1}---\autoref{fig3} give the comparison
between the target shape with iterated shape for the viscosity
coefficients $\alpha=0.1, 0.01$ and $0.001$, respectively.

\begin{figure}[!htbp]
\centering
  \includegraphics[width=.8\textwidth]{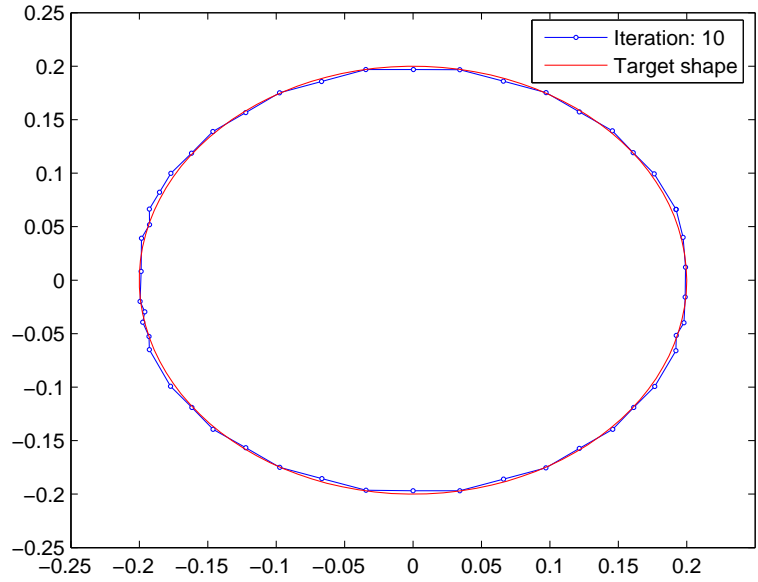}
  \caption{Case 1: $\alpha=0.1$, CPU time: 111.984 s.\label{fig1}}
\end{figure}

\begin{figure}[!htbp]
\centering
  \includegraphics[width=.8\textwidth]{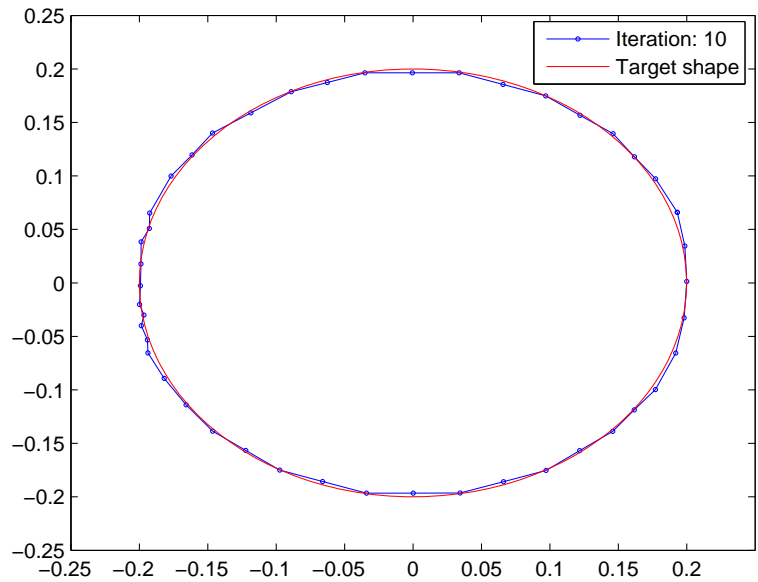}
  \caption{Case 1: $\alpha=0.01$, CPU time: 125.969 s.\label{fig2}}
\end{figure}

\begin{figure}[!htbp]
\centering
  \includegraphics[width=.8\textwidth]{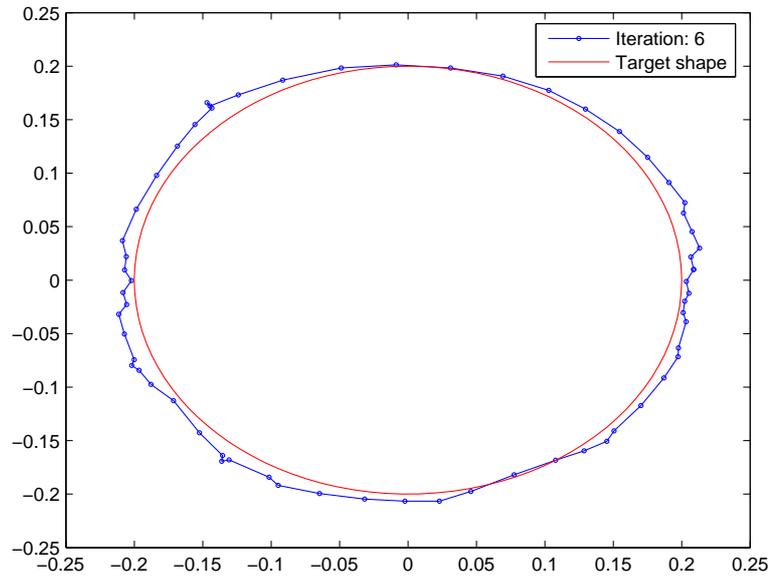}
  \caption{Case 1: $\alpha=0.001$ with $h=20$, CPU time: 185.578 s.\label{fig}}
\end{figure}

\begin{figure}[!htbp]
\centering
  \includegraphics[width=.8\textwidth]{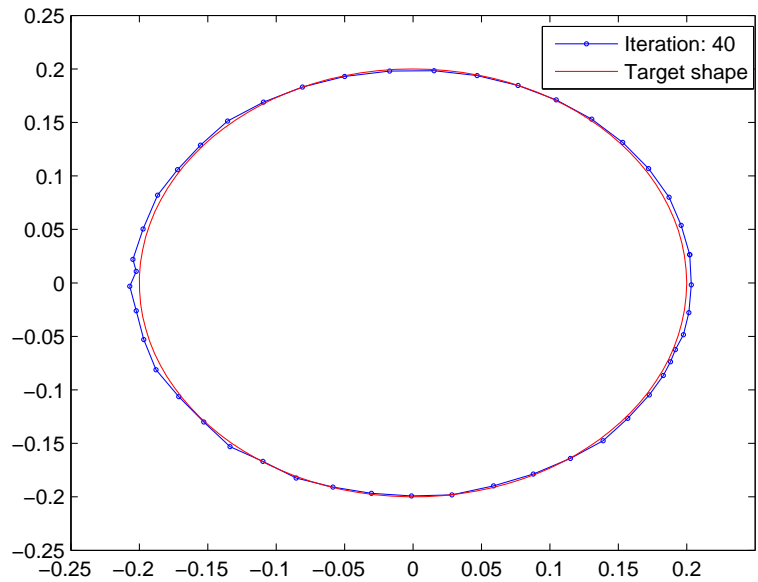}
  \caption{Case 1: $\alpha=0.001$ with $h=5$, CPU time: 700.016 s.\label{fig3}}
\end{figure}

For $\alpha=0.1$ and $\alpha=0.01$, we choose the initial step
$h=20$ and in \autoref{fig1} and \autoref{fig2}, we give the final
shape at iteration 10 with CPU times. However for $\alpha=0.001$,
one can not obtain a good result when $h=20$ (see \autoref{fig}).
Thus we should reduce the initial descent step, and then in
\autoref{fig3}, we give the final shape at iteration 40 with the
initial step $h=5$. By comparison with \autoref{fig}, we find that
though we need much more CPU time for $h=5$, but it have a nicer
reconstruction.

\autoref{fig4} represents the fast convergence of the cost
functional for the various viscosities in Case 1.
\begin{figure}[!htbp]
\centering
  \includegraphics[width=.89\textwidth]{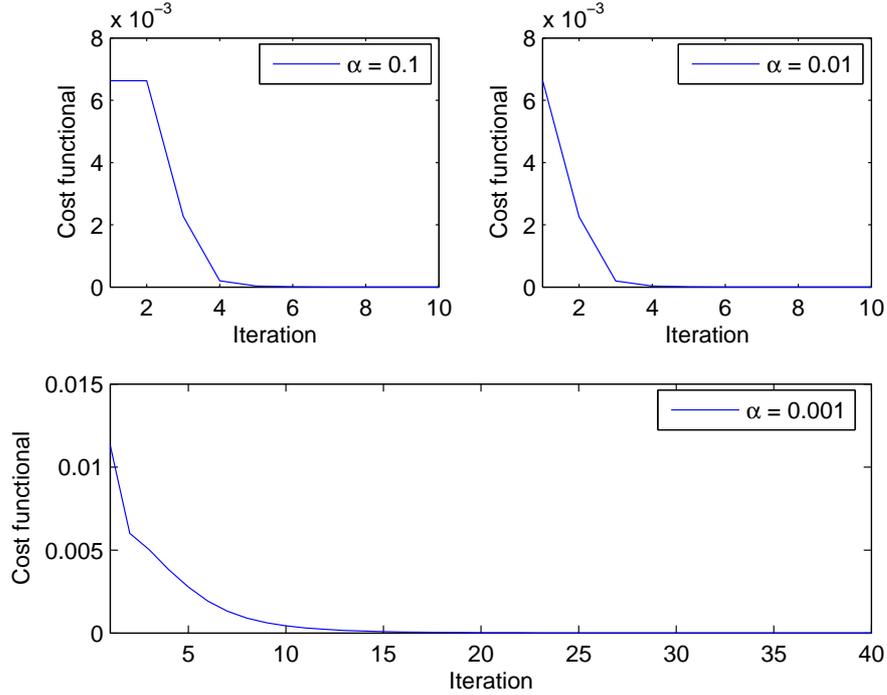}
  \caption{Convergence history in Case 1 for $\alpha=0.1, 0.01$ and $0.001$.\label{fig4}}
\end{figure}
\begin{figure}[!htbp]
\centering
  \includegraphics[width=.8\textwidth]{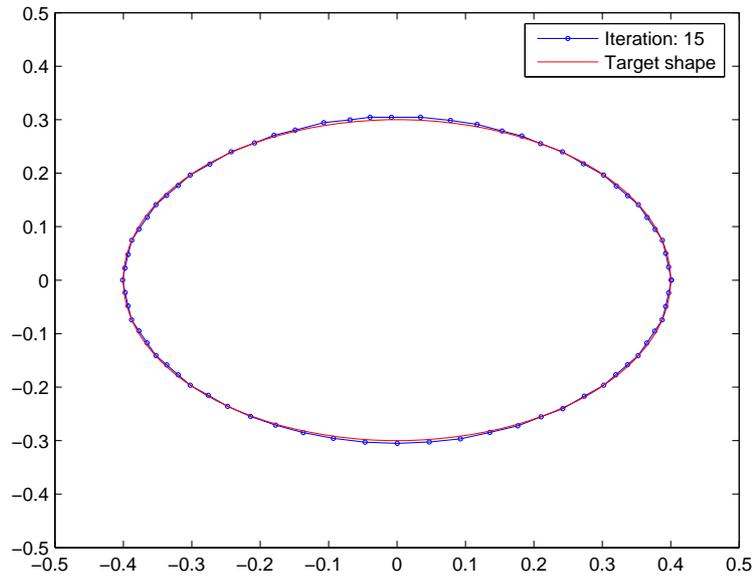}
  \caption{Case 2: $\alpha=0.1$, CPU time: 196.609 s.\label{fig5}}
\end{figure}

\begin{figure}[!htbp]
\centering
  \includegraphics[width=.8\textwidth]{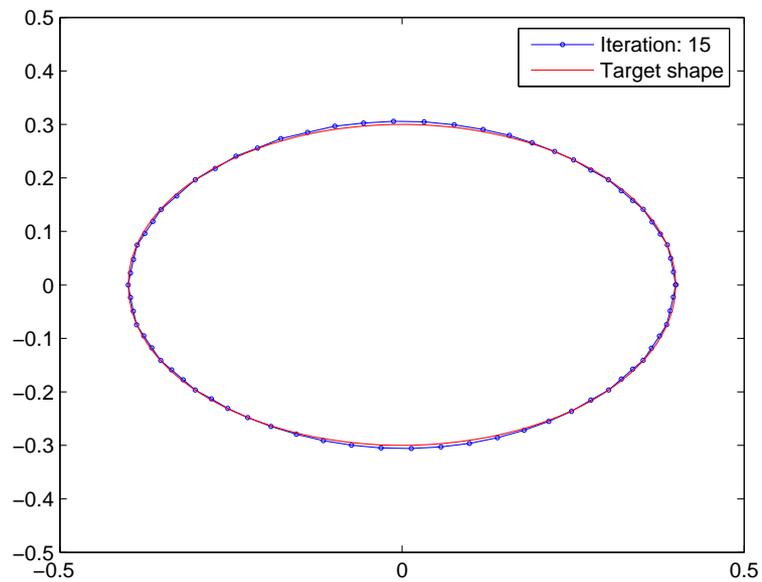}
  \caption{Case 2: $\alpha=0.01$, CPU time:  207.469 s.\label{fig6}}
\end{figure}

\begin{figure}[!htbp]
\centering
  \includegraphics[width=.8\textwidth]{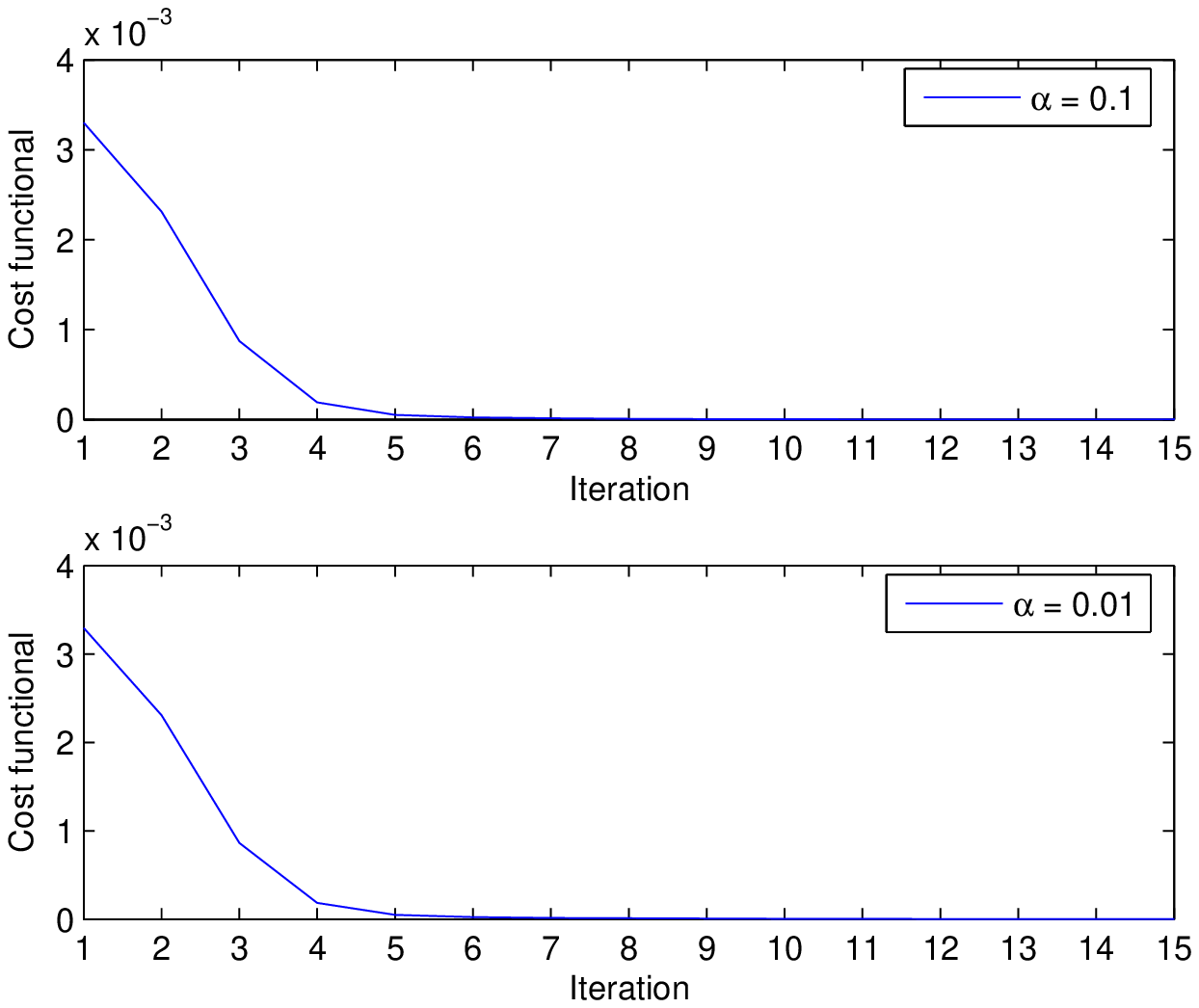}
  \caption{Convergence history in Case 2 for $\alpha=0.1, 0.01$.\label{fig7}}
\end{figure}


 In Case 2,
\autoref{fig5} and \autoref{fig6} show the comparison between the
target shape with iterated shape at iteration 15 for the viscosity
$\alpha=0.1$ and $\alpha=0.01$, respectively. At this time, we
choose the initial step $h=20$.  We also give the CPU run times at
the 15 iterations for $\alpha=0.1$ and $\alpha=0.01$. Unfortunately,
we can not get a good reconstruction for $\alpha=0.001$ in this
case. \autoref{fig7} gives the convergence history of the cost
functional $J(\oo)$ for $\alpha=0.1, 0.01$.

Finally, we can conclude that the proposed gradient type algorithm
is an efficient one in both of our test cases. Unfortunately for
large Reynold numbers, we can not obtain the nice results quickly.
Hence further research is necessary on efficient implementations for
very large Reynold numbers and real problems in the industry.

\end{document}